\documentclass[11pt]{article} 

\usepackage{amssymb}
\usepackage{amsmath}
\usepackage{CJKutf8} 
\usepackage{geometry}
\usepackage{braket}
\usepackage{pdflscape}
\usepackage{subfig}
\usepackage{url}
\usepackage{algorithm}
\usepackage{algpseudocode}
\usepackage{setspace}
\usepackage{appendix}
\usepackage{framed}
\usepackage{xcolor}
\usepackage{graphicx} 
\usepackage{multirow}
\usepackage{amsthm}
\usepackage{graphicx}
\usepackage{hyperref}
\usepackage{booktabs}
\usepackage{tabularx}
\usepackage[utf8]{inputenc}
\usepackage{nicematrix}
\newtheorem{theorem}{Theorem}
\newtheorem{lemma}[theorem]{Lemma}
\newtheorem{example}{Example}

\newtheorem{corollary}[theorem]{Corollary}
\newtheorem{assumption}{Assumption}

\newcommand{\R}{\mathbb{R}}

\newcommand{\Rnn}{\mathbb{R}^{n \times n}}
\newcommand{\Rmn}{\mathbb{R}^{m \times n}}

\newcommand{\Rnr}{\mathbb{R}^{n \times r}}

\newcommand{\Rmr}{\mathbb{R}^{m \times r}}

\newcommand{\norm}[1]{\left\| #1 \right\|}

\renewcommand\arraystretch{1.2}

\newcommand{\M}{\mathcal M}
\newcommand{\Mr}{\mathcal{M}_r}

\newcommand{\Tr}{\mathcal{T}_r}

\newcommand{\proj}[2]{\mathcal{P}_{#1} \left[ #2 \right]}

\newcommand{\PPDEIM}{\mathcal P^{\angle}}

\usepackage{accents}
\newcommand*{\dt}[1]{%
  \accentset{\mbox{\large\bfseries .}}{#1}}

\newcommand{\reflectangle}{\mathbin{\mathpalette\reflectangleaux\angle}}
\newcommand{\reflectangleaux}[2]{\reflectbox{$#1#2$}}

\usepackage{authblk}
\renewcommand\Affilfont{\small\normalfont}   
\makeatletter
\renewcommand\AB@affilsepx{\\\protect\Affilfont}  
\makeatother

\newcommand{\msc}[1]{\par\smallskip\noindent\textbf{MSC 2020:} #1\par}
\newcommand{\kw}[1]{\noindent\textbf{Keywords:} #1\par}

\usepackage{fancyhdr}
\title{Interpolatory Dynamical Low-Rank Approximation: \\
Theoretical Foundations and Algorithms}

\author[1]{Benjamin Carrel}
\author[2]{Daniel Kressner}
\author[2]{Hei Yin Lam}
\author[3]{Bart Vandereycken}

\affil[1]{Laboratory for Simulation and Modelling, PSI, Villigen, Switzerland\newline
\texttt{benjamin.carrel@psi.ch}}
\affil[2]{Institute of Mathematics, EPFL, Lausanne, Switzerland\newline 
\texttt{\{hysan.lam, daniel.kressner\}@epfl.ch}}
\affil[3]{Section of Mathematics, University of Geneva, Geneva, Switzerland\newline 
\texttt{bart.vandereycken@unige.ch}}

\date{\today}

\begin{document}
\maketitle

\begin{abstract}
Dynamical low-rank approximation (DLRA) is a widely used paradigm for solving large-scale matrix differential equations, as they arise, for example, from the discretization of time-dependent partial differential equations on tensorized domains. Through orthogonally projecting the dynamics onto the tangent space of a low-dimensional manifold, DLRA achieves a significant reduction of the storage required to represent the solution. However, the need for evaluating the velocity field can make it challenging to attain a corresponding reduction of computational cost in the presence of nonlinearities. In this work, we address this challenge by
replacing orthogonal tangent space projections with oblique, data-sparse projections selected by a discrete empirical interpolation method (DEIM).
At the continuous-time level, this leads to DLRA-DEIM, a well-posed differential inclusion (in the Filippov sense) that captures the discontinuities induced by changes in the indices selected by DEIM. 
We establish an existence result, exactness property and error bound for DLRA-DEIM that match existing results for DLRA.
 For the particular case of QDEIM, a popular variant of DEIM, we provide an explicit convex-polytope characterization of the differential inclusion.
Building on DLRA-DEIM, we propose a new class of 
projected integrators, called PRK-DEIM, that combines explicit
Runge--Kutta methods with DEIM-based projections.
We analyze the convergence order of PRK-DEIM and show that it matches the accuracy of previously proposed projected Runge-Kutta methods, while being significantly cheaper. Extensions to exponential Runge--Kutta methods and low-order tensor differential equations demonstrate the versatility of our framework.
\end{abstract}

\msc{65F55, 65L20, 65L05, 15A69, 53Z50}
\kw{dynamical low-rank approximation; DEIM; projected Runge–Kutta; tensor methods}

\section{Introduction} \label{sec: introduction}

This work is concerned with the numerical integration of large-scale matrix differential equations of the form
\begin{equation} \label{eq: full order model}
  \dt{A}(t) = F(A(t)), \quad A(0) = A_0 \in \Rmn.
\end{equation}
Such differential equations arise, for example, from spatial discretizations of time-dependent partial differential equations~\cite{leveque2007finite,quarteroni1994numerical,trefethen2000spectral} on tensorized domains and in control theory \cite{engwerda2005lq,lewis2012optimal}. 
The full-order model~\eqref{eq: full order model} becomes computationally expensive to solve when $m$ and $n$ are large. 
Model order reduction techniques aim at reducing computational cost and memory requirements by exploiting some form of data compression. A common situation is that
$A(t)$ admits low-rank data compression, because of, e.g., spatial smoothness.
Dynamical low-rank approximation (DLRA) is built on this observation and has demonstrated its versatility for a broad variety of large-scale problems in, e.g., plasma physics \cite{cassini2022efficient,coughlin2022efficient,einkemmer2024accelerating,einkemmer2018low}, radiation therapy \cite{kusch2023robust}, chemical kinetics \cite{einkemmer2024low,jahnke2008dynamical,prugger2023dynamical}, wave propagation~\cite{hochbruck2023rank,zhao2025low}, uncertainty quantification \cite{donello2023oblique,feppon2018dynamically,kazashi2025dynamical,nobile2025petrov}, and machine learning \cite{schotthofer2024federated,schotthofer2025dynamic,schotthofer2022low}.

DLRA evolves an approximate solution to \eqref{eq: full order model} on the smooth embedded sub-manifold
$$\Mr = \{X \in \Rmn \mid \mathrm{rank}(X) = r \}$$
for some fixed rank $r \ll m,n$.
This evolution is achieved by first compressing the initial matrix $A_0$ to a rank-$r$ matrix $X_0$ and then using tangent space projection to ensure that the dynamics remains
on $\Mr$:
\begin{equation} \label{eq: DLRA}
  \dt{X}(t) = \proj{X(t)}{F(X(t))}, \quad X(0) = X_0 \in \Mr.
\end{equation}
Here, $\proj{X}{\,\cdot\,} \colon \Rmn \rightarrow T_X \Mr$ is the orthogonal projection onto the tangent space 
$T_X \Mr$ at $X \in \Mr$. 
Using that any such $X$ admits a rank-$r$ factorization of the form $X = USV^\top$
with $S \in \R^{r\times r}$ and 
orthonormal\footnote{Throughout this work we say that a rectangular tall matrix is orthonormal if its columns are orthonormal.} matrices $U \in \Rmr$, $V \in \Rnr$, one can write
\begin{equation} \label{eq: orthogonal tangent projection}
\proj{X}{Z} = UU^\top Z + Z VV^\top - UU^\top Z VV^\top.
\end{equation}

The expression~\eqref{eq: orthogonal tangent projection} allows one to rewrite~\eqref{eq: DLRA} as a (low-dimensional) system of differential equations for the rank-$r$ factors and, in principle, one can simply apply a standard numerical integrator to it. However, such an approach has been observed to lead to numerical instabilities in the presence of small singular values. To address this issue, specific methods~\cite{lubich2014projector,ceruti2024robust,ceruti2022unconventional,kieri2019projection,kusch2024second} have been developed over the last decade, and this remains an active area of research.

Another challenge of DLRA is that the memory reduction achieved by representing $X$ in terms of its rank-$r$ factors does not necessarily translate into a corresponding reduction of computational cost. In particular, solving~\eqref{eq: DLRA} involves the repeated orthogonal tangent projection of $F(X)$. Unless $F$ is linear and/or highly structured, this requires
evaluating all entries of $F$ at the explicitly formed matrix $X$. Even a highly structured nonlinearity can incur significant cost. Take, for example, the cubic nonlinearity
$X * X * X$, where $*$ denotes the elementwise product. If $\mathrm{rank}(X) = r$ then $\mathrm{rank}(X * X * X) \leq r^3$, indicating a massive intermediate rank growth that may annihilate all computational gains expected from DLRA. In this work, we will explore the use of DEIM (discrete empirical interpolation method) for addressing this challenge.

DEIM~\cite{barrault2004empirical,chaturantabut2010nonlinear} is a popular technique for efficiently evaluating nonlinear terms in the context of model order reduction methods such as POD (proper orthogonal decomposition). To recall the basic principles of DEIM, let us consider a (vector) ODE of the form 
$$\dt{x}(t) = f(x(t)), \qquad x(0) = x_0 \in \R^n,$$
with $f \colon \R^n \rightarrow \R^n$.
Suppose one has an orthonormal basis $U = [u_1, \ldots, u_r] \in \Rnr$ with $r\ll n$ (obtained from, e.g., previous simulations) that captures most of the dynamics very well.
POD then applies the orthogonal projection $UU^\top$ to both sides of the ODE, resulting in the approximation 
$x(t) \approx U y(t)$, where 
$$\dt{y}(t) = U^T f(U y(t)), \qquad y(0) = U^\top x_0 \in \R^r.$$
To reduce the cost of evaluating $f$, DEIM performs an oblique instead of an orthogonal projection by choosing $r$ rows of $U$. The corresponding row indices $p_1, \ldots, p_r \in \{1,\ldots, n\}$ define 
the selection matrix
$$S_U = [\pmb e_{p_1}, \ldots, \pmb e_{p_r}] \in \R^{n \times r},$$
where $\pmb e_i \in \R^n$ denotes the $i$th unit vector. Applying the oblique projection \[P^{\angle}_U:=U ( S_U^\top U)^{-1} S_U^\top\] to the ODE results in the approximation
$x(t) \approx U z(t)$ with
$$\dt{z}(t) = (S_U^\top U)^{-1} S_U^\top f(Uz(t)), \qquad z(0) = (S_U^\top U)^{-1} S_U^\top x_0 \in \R^r.$$
The key advantage of this POD-DEIM approximation is that only $r$ entries of $f$ need to be evaluated in order to form $S_U^\top f$.
The quality of this approximation critically depends on a good choice of indices.

\subsection*{Contributions and outline}

There have been efforts to bring the benefits of DEIM to DLRA. Notably, the TDB-CUR method from~\cite{donello2023oblique} applies DEIM procedures to the columns and rows of the full-order model~\eqref{eq: full order model} discretized in time by, e.g., an explicit Runge–Kutta method. Avoiding the disadvantages mentioned above, TDB-CUR
can be much faster than DLRA on problems containing non-linear terms since it only needs to evaluate the vector field at certain rows and columns, selected by DEIM. More closely related to the current work is the interpolatory projector-splitting integrator from~\cite{dektor2025collocation, dektor2024interpolatory} where DEIM is used to replace the orthogonal projectors in the projector-splitting integrators from~\cite{lubich2014projector, lubich2015time_OV}. Also for these methods, large computational savings were shown for vector fields with non-linear terms. From a theoretical perspective, the effects of such DEIM procedures on the dynamics and the time stepping error are not fully understood. In particular, there is no model that meaningfully captures DEIM as the time-step size converges to $0$. In addition, there is no stability or convergence analysis available for the discrete versions of these DLRA-with-DEIM integrators.

A key contribution of this work is to bridge this gap; we consider a continuous-time model called DLRA-DEIM that incorporates DEIM into DLRA. Replacing orthogonal by oblique projections in~\eqref{eq: DLRA}--\eqref{eq: orthogonal tangent projection} leads to
\begin{equation} \label{eq:dlardeimintro}
\dt{X}(t) = P_{U(t)}^{\angle} F(X(t)) - P_{U(t)}^{\angle} F(X(t)) P_{V(t)}^{\reflectangle} + F(X(t)) P_{V(t)}^{\reflectangle}, 
\end{equation}
with $P^{\reflectangle}_{V(t)}:=(P^{\angle}_{V(t)})^\top$ and the factorized representation $X(t) = U(t) S(t) V(t)^\top$.
However, a fundamental difficulty is that the re-adjustment of selected indices in DEIM-like procedures introduces time discontinuities, which make the differential equation~\eqref{eq:dlardeimintro} ill-posed. In Section~\ref{sec:theory}, 
we explain how one can use differential inclusions to make sense of~\eqref{eq:dlardeimintro} and provide a detailed analysis when QDEIM is used for index selection. Inspired by~\eqref{eq:dlardeimintro}, Section~\ref{sec: PRK-DEIM} develops and analyzes PRK-DEIM, which incorporates DEIM into the projected Runge-Kutta methods from~\cite{kieri2019projection}. Our numerical experiments confirm that PRK-DEIM can reduce execution times significantly, while only having a negligible impact on the error. Finally, in Section~\ref{sec:tucker}, we extend PRK-DEIM from matrix to (low-order) tensor differential equations, using a combination of the Tucker format and DEIM. Compared to other DLRA methods with DEIM, like~\cite{donello2023oblique}, our method is less intrusive since parts of the vector field that map onto the tangent space can be treated exactly without interpolation by DEIM. In case of the discretized Laplacian, for example, such interpolation is relatively costly due to the dependence of the sampled rows/columns on non-sampled rows/columns.

\section{Preliminaries}

\subsection{Assumptions for the dynamical low-rank approximation}

In the following, we recall assumptions commonly made in the analysis of ODEs and DLRA, which are also needed in our analysis. The following classical assumption guarantees the existence and uniqueness of the solution to the full-order model~\eqref{eq: full order model}.
\begin{assumption} \label{ass: Lipschitz}
The field $F: \Rmn \rightarrow \Rmn$ is Lipschitz continuous, that is, there exists a constant $L \ge 0$ such that 
$\norm{F(X) - F(Y)}_F \leq L \norm{X - Y}_F$ holds for all $X, Y \in \Rmn$.
\end{assumption}

The concept of one-sided Lipschitz continuity can yield tighter error bounds, especially when $\ell < 0$.
\begin{assumption}[One-sided Lipschitz constant]\label{ass: one-sided F}
The function $F$ is one-sided Lipschitz continuous with constant $\ell \in \mathbb  R$,
that is,
$\langle X-Y,F(X)-F(Y) \rangle \leq \ell \norm{X-Y}_F^2$ holds for all $X, Y \in \Rmn$.
\end{assumption}

The analysis of DLRA usually relies on the following assumption from \cite{koch2007dynamical}.
\begin{assumption}[Small normal component] \label{ass: epsilon close}
The normal components of $F(Y)$ are $\varepsilon_r$-small, that is,
    $\|\mathcal{P}_{Y}[F(Y)]- F(Y)\|_F\leq \varepsilon_r$ holds
    for all $Y\in \mathcal{M}_r$.
\end{assumption}
\begin{theorem}[Accuracy of DLRA~\cite{kieri2019projection}]
\label{thm: DLRA_acc}
Under Assumptions \ref{ass: one-sided F} and \ref{ass: epsilon close}, the difference between the solution $X(t)$ to the DLRA~\eqref{eq: DLRA} and the solution $A(t)$ to the full-order model~\eqref{eq: full order model} satisfies
\begin{equation} \label{eq: error made by DLRA}
\norm{A(t) - X(t)}_F \leq e^{t \ell} \norm{A_0 - X_0}_F + \varepsilon_r \cdot \int_0^t e^{(t-s) \ell}\,\mathrm{d}s.
\end{equation}
\end{theorem}
The analysis of the projected Runge-Kutta methods from~\cite{kieri2016discretized}  makes the following additional assumption.
\begin{assumption}[Bounded vector field] \label{ass: bounded F}
There exists a constant $B \in \R$ such that $\norm{F(X)}_F \leq B$ for all $X \in \Rmn$.
\end{assumption}
For convenience and to avoid technicalities, the assumptions above have a global nature, that is, they are supposed to hold everywhere in $\Rmn$ or in $\mathcal{M}_r$. As explained in~\cite[p.1025]{kieri2016discretized}, these assumptions can be weakened to local conditions that only hold in some suitably large tubular neighborhood of the exact solution $A(t)$ for $0 \leq t \leq T$.

\subsection{DEIM procedures}

As explained in the introduction, our construction of DLRA-DEIM centrally involves oblique projectors of the form $U (S_U^\top U)^{-1} S_U^\top$ for a selection matrix $S_U$ and an orthonormal matrix $U$. Given an arbitrary matrix $Z\in \R^{m\times n}$, one has~\cite{chaturantabut2010nonlinear} that
\begin{equation} \label{eq:quasioblique}
  \|Z - P_U^{\angle} Z \|_F \le \| (S_U^\top U)^{-1} \|_2 \| Z - UU^\top Z \|_F.
\end{equation}
Hence, the oblique projection of the columns of $Z$ onto the range of $U$ leads to an approximation error that differs from the best approximation error by a factor at most $\| (S_U^\top U)^{-1} \|_2$. This factor determines the quality of the indices $p_1,\ldots, p_r$ in the selection matrix $S_U$. Table~\ref{tab:algo_comparison} provides a brief overview of popular DEIM procedures and their known upper bounds 
\begin{equation}\label{eq:def Cmr}
\| (S_U^\top U)^{-1} \|_2 \le C_{m,r}.
\end{equation}
For example, the original DEIM~\cite{chaturantabut2010nonlinear} is a greedy selection strategy with an upper bound that grows exponentially with $r$. Let us stress that this exponential growth is extremely rare and has not been observed to cause severe issues in practice.

\begin{table}[h]
  \centering
  \renewcommand{\arraystretch}{1.3} 
  \begin{tabularx}{\textwidth}{l >{\centering\arraybackslash}p{3.5cm} X}
    \toprule
    \textbf{DEIM procedure} & \textbf{$C_{m,r}$} & \textbf{Remarks} \\ 
    \midrule
    original DEIM~\cite{chaturantabut2010nonlinear} & 
    $\left( 1+ \sqrt{2m} \right)^{r-1} \|u_1\|_{\infty}^{-1}$ & Original greedy algorithm. \\
    
    QDEIM \cite{drmac2016new} & 
    $\sqrt{m - r - 1} \frac{\sqrt{4^r + 6r - 1}}{3}$ & Pivoted QR decomposition of $U^\top$. \\
    
    SRRQR \cite{drmac2018discrete,gu1996efficient, osinsky2018rectangular} & 
    $\sqrt{1 + \eta^2 r (m-r)}$ & Enhances QDEIM by strong rank-revealing QR (SRRQR) with parameter $\eta > 1$ of $U^\top$. Worst case cost is $O\left( m r^3 \log_{\eta}(r)\right)$ operations but $O(mr^2)$ in practice with moderate $\eta$. \\
    
    Osinsky~\cite{osinsky_close_2023} &  
    $\sqrt{1 + r(m-r)}$ & Does not leverage existing implementations of pivoted QR. \\
    
    {ARP} \cite{cortinovis_adaptive_2024} &  
    $\sqrt{1 + r(m-r)}$ \, \, \emph{(in expectation)} & 
    Randomized version of Osinsky, based on adaptive leverage score sampling.  \\ 
    
    Gappy POD+E \cite{peherstorfer2020stability} & 
    N.A. & Requires $O(m r^2 + \ell (r+\ell)^3)$ operations, with oversampling parameter $\ell$  for enhanced stability.\\
    \bottomrule
  \end{tabularx}
  \caption{Overview of DEIM procedures for selecting $r$ rows from an orthonormal matrix $U=[u_1,\ldots, u_r] \in \mathbb{R}^{m \times r}$, along with known upper bounds $C_{m,r}$ on $\|(S_U^\top U)^{-1}\|_2$. Unless otherwise mentioned, the algorithms require $O(m r^2)$ operations. }
  \label{tab:algo_comparison}
\end{table}

\begin{figure}[ht!]
\centering
\includegraphics[width=0.9\textwidth]{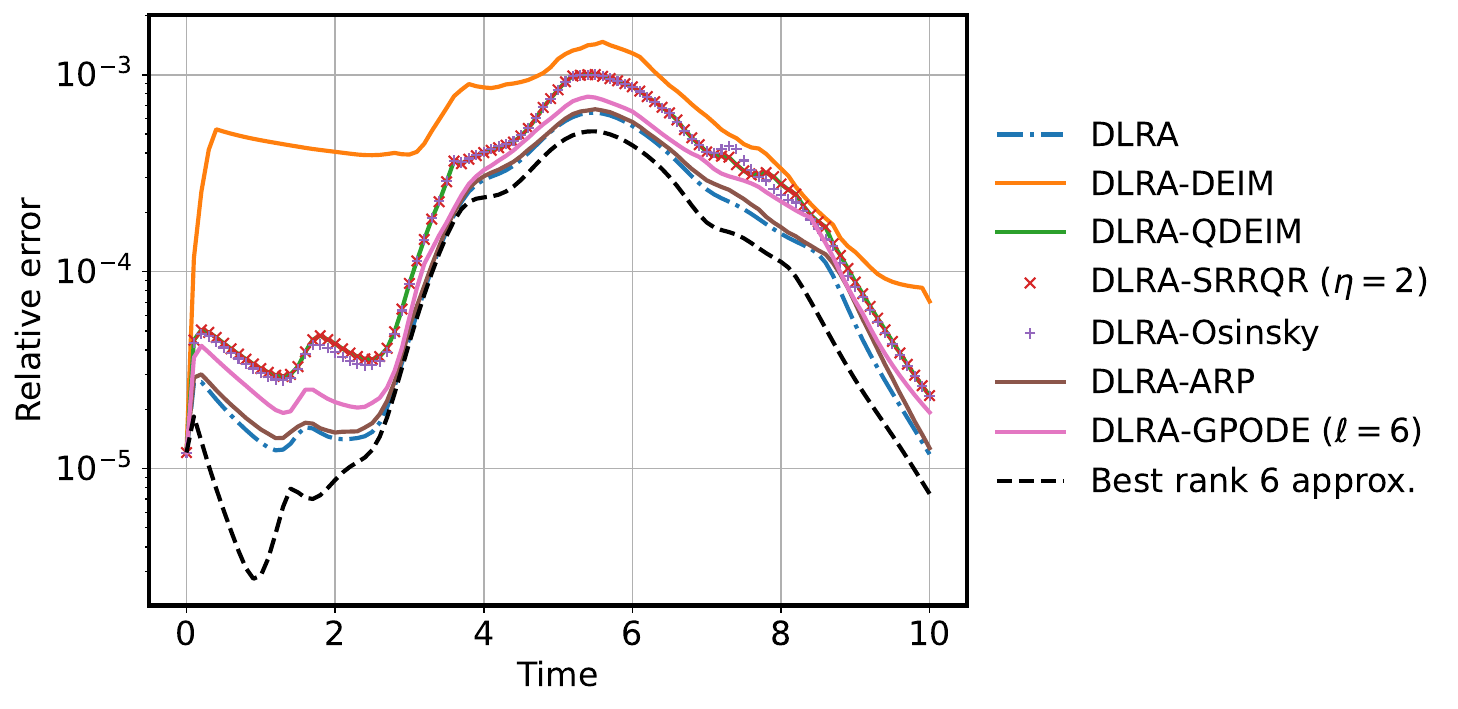}
\caption{Relative error between reference solution and DLRA, using the different DEIM techniques from Table \ref{tab:algo_comparison}, for a discretized Allen-Cahn equation.
\label{fig: motivating figure}}
\end{figure}

In a preliminary numerical study, we have tested all DEIM procedures of Table~\ref{tab:algo_comparison} in combination with the numerical integrators for DLRA-DEIM discussed in Section~\ref{sec: PRK-DEIM}. For generating the error plots of Figure~\ref{fig: motivating figure}, we have applied a second order projected Runge--Kutta method
to the Allen-Cahn equation \eqref{eq: Allen-Cahn} with $\kappa = 0.01$ and discretized problem size $n=64$. The time step size $h = 10^{-3}$ is sufficiently small so that the time integration error does not visibly affect the overall error.  DLRA corresponds to using orthogonal tangent space projections~\eqref{eq: orthogonal tangent projection} with rank $r=6$ and one observes that its error remains quite close to the one of the best rank-$6$ approximation to the reference solution of the full-order model.
Using oblique instead of orthogonal projections introduces additional error; while not drastic, the differences between the different DEIM procedures can be significant---up to an order of magnitude. Based on these findings, our subsequent numerical experiments will focus on SRRQR and ARP.

\subsection{A closer look at QDEIM}

The discrete choice of indices in DEIM introduces discontinuities. In the following, we will provide a detailed characterization of these discontinuities for QDEIM, which is a popular choice and amenable to such an analysis. We suspect that the analysis is similar for the original DEIM, SRRQR and Osinsky's procedures, while randomization and oversampling add significant complications for the analysis of the other DEIM algorithms.

Algorithm~\ref{alg:QDEIM} recalls QDEIM from~\cite{drmac2016new} for the more general case of a full-rank matrix $U$. We use Matlab's colon notation; $U^{k-1}(j,:)$ denotes the $j$th row of the matrix $U^{k-1}$.  QDEIM chooses in each step the row of maximum norm and then removes this row by an orthogonal projection. Note that this is equivalent to applying QR with column pivoting~\cite[Sec 5.4.2]{GolubVanLoan2013} to $U^\top$, which yields a simple and numerically stable way of carrying out QDEIM; see~\cite[Section 2.1.1]{drmac2016new} for more details.

There may be several rows having maximal norm in any of the $r$ steps of QDEIM. To break such ties, we select the row of maximal norm having the smallest row index $p_k$. 
\begin{algorithm}[H]
\caption{QDEIM with tie-breaking rule}
\label{alg:QDEIM}
\begin{algorithmic}[1]
\Require Matrix $U\in\Rmr$ of full column rank. 
\State Set $U^0=U$.
\For {$k=1,\ldots, r$}
\State Find smallest index $p_k$ such that $\|U^{k-1}(p_k,:)\|_2=\max_{j\in \{1,\ldots m\}}\|U^{k-1}(j,:)\|_2.$ \label{line: max}
\State Set $u={U^{k-1}(p_k,:)^\top}/{\|U^{k-1}(p_k,:)\|_2}.$
\State Update $U^k\leftarrow U^{k-1}\big(I-uu^\top\big)$. \label{line: update}
\EndFor
\Ensure Selection matrix $S^{\mathsf{Q}}_U = [\pmb e_{p_1}, \ldots, \pmb e_{p_r}]$.
\end{algorithmic}
\end{algorithm}

\begin{example} \label{ex:small}
Consider the following orthonormal matrix $U(t)$ depending continuously on $t \in [0, 1]$:
$$U(t) =  \begin{bmatrix} t &  0 \\ \sqrt{1-t^2} & 0 \\ 0 & 1 \end{bmatrix}.$$
For every sufficiently small $h>0$, Algorithm~\ref{alg:QDEIM} returns
        \begin{equation}
        \label{eq:selection matrix}
            S^{\mathsf{Q}}_{U(1/\sqrt{2} - h)} = \begin{bmatrix} 0 & 0 \\ 0 & 1 \\ 1 & 0 \end{bmatrix}, \quad S^{\mathsf{Q}}_{U(1/\sqrt{2} + h)} = \begin{bmatrix} 0 & 1 \\ 0 & 0 \\ 1 & 0 \end{bmatrix}.\end{equation}
  For $h=0$, Algorithm~\ref{alg:QDEIM} encounters a tie in the second loop, which is resolved in Line~\ref{line: max} by choosing the smaller index $p_2 = 1$,  leading to $S^{\mathsf{Q}}_{U(1/\sqrt{2})}=S^{\mathsf{Q}}_{U(1/\sqrt{2}+h)}$. QDEIM thus exhibits a discontinuity at $t=1/\sqrt{2}$, and so does its corresponding oblique projection.        
\end{example}

The example above conveys the intuition that any tie in QDEIM leads to a discontinuity. The following theorem formalizes this statement.
\begin{theorem} \label{thm:QDEIM}
Consider the function $S^{\mathsf{Q}}_{(\cdot)}: \R^{m\times r}\to \R^{m\times r}$ that maps a full-rank matrix $U$ to the selection matrix $S^{\mathsf{Q}}_U$ returned by QDEIM, Algorithm~\ref{alg:QDEIM}.
Then $S^{\mathsf{Q}}_{(\cdot)}$ is continuous at $U$ if and only if Line~\ref{line: max} of Algorithm~\ref{alg:QDEIM} does not encounter a tie when selecting the $r$ indices. 
\end{theorem}
Because selection matrices form a discrete set, $S^{\mathsf{Q}}_U$ remains constant in a sufficiently small neighborhood of $U$ when the continuity condition of Theorem~\ref{thm:QDEIM} holds.

\begin{proof}
We first note that the full-rank condition of Algorithm~\ref{alg:QDEIM} remains satisfied in a sufficiently small neighborhood of $U$.

To prove the implication $(\Leftarrow)$, suppose that Algorithm~\ref{alg:QDEIM} applied to $U$ encounters no ties in Line~\ref{line: max}. 
We proceed by induction and assume that $p_1, \ldots ,p_{k-1}$ and $U^1, \ldots , U^{k-1}$ are continuous in $U$. Consider the $k$th iteration of Algorithm~\ref{alg:QDEIM}. Because there is no tie in Line~\ref{line: max} and the Euclidean norm is continuous, the pivot $p_k$ will remain the same (and hence continuous) in a small neighborhood of $U$. Since $U$ has full rank, the row $U^{k-1}(p_k,:)$ is non-zero and, in turn, the update in Line~\ref{line: update} of $U^k$ is also continuous. By induction, this shows the statement for $k=r$ and we have proved that $S^{\mathsf{Q}}_{(\cdot)}$ is continuous at $U$.

To prove $(\Rightarrow)$, we show the contrapositive. Suppose that Algorithm~\ref{alg:QDEIM} encounters ties and the first tie appears in iteration $k$. Then there exist $i,j\in \{1,\ldots, m\}\setminus\{p_1, \ldots,p_{k-1}\}$ with $i\neq j$ such that 
$$\|U^{k-1}(i,:)\|_2=\|U^{k-1}(j,:)\|_2=\max_{\ell \in \{1,\ldots, m\}}\|U^{k-1}(\ell,:)\|_2.$$ Now, consider a perturbation $U'=U+ \epsilon \cdot \pmb{e}_{i}U^{k-1}(i,:)$ for some $\epsilon\in \R$ and  let $U'^{\ell}$, $p'_\ell$ denote the result of applying $\ell$ iterations of Algorithm~\ref{alg:QDEIM} to $U'$. From the continuity established in the first part of the proof,
it follows that the rows selected until iteration $k-1$ remain unaltered for sufficiently small $\epsilon$, that is,
$$p'_{\ell} = p_\ell, \quad U'^{\ell-1}(p'_{\ell},:)=U^{\ell-1}(p_{\ell},:)\quad \text{for}\quad \ell=1, \ldots,k-1.$$ By the update in Line~\ref{line: update}, $U^{k-1}(i,:)$ is orthogonal to $U^0(p_1,:),U^1(p_2,:)$, $\ldots$, $U^{k-2}(p_{k-1},:)$ and thus $U'^{k-1}=U^{k-1}+ \epsilon \cdot \pmb{e}_{i}U^{k-1}(i,:)$. We now have  $\|U'^{k-1}(i,:)\|_2=(1+ \epsilon)\|U^{k-1}(i,:)\|_2$. Then, $p_k=i$ when $\epsilon>0$ and $p_k\neq i$ when $\epsilon<0$. Hence, the pivot selection is discontinuous at $U$.
\end{proof}

Our analysis of DLRA-DEIM will require the following characterization of all possible selections in the vicinity of a matrix $U$ (restricted to the Stiefel manifold of orthonormal matrices):
\begin{equation} \label{eq:vicinity}
\mathcal{S}_U := \big\{
S \;\big|\; 
\exists \{U_i \in \R^{m\times r}\} \text{ such that } \, U^\top_i U_i = I_r,\, U_i \to U, \text{ and }\, S^{\mathsf{Q}}_{U_i} \to S
\big\}.
\end{equation}
For Example~\ref{ex:small}, it is not difficult to see that $\mathcal{S}_{U(1/\sqrt{2})}$ consists of the two selection matrices shown in~\eqref{eq:selection matrix}. On the other hand, we can also consider the effect of switching off tie breaking in Algorithm~\ref{alg:QDEIM}, that is, Line~\ref{line: max} chooses any arbitrary row of maximal norm:
\begin{equation} \label{eq:rule}
\|U^{k-1}(p_k, :)\|_2 = \max_{j \in \{1, \ldots, m\}} \|U^{k-1}(j, :)\|_2.
\end{equation}
Note that the subsequently updated matrices $U^\ell$, $\ell \ge k$, and indices are not invariant under this choice. The following set collects \emph{every} selection matrix that results from a feasible path of indices $p_1, p_2, \ldots, p_r$ (each satisfying~\eqref{eq:rule}):
\begin{equation} \label{eq:allselections}
\mathcal{I}_U = \left\{ S \;\middle|\; S \text{ is a feasible selection matrix of QDEIM without tie-breaking applied to } U  \right\}.
\end{equation}
For Example~\ref{ex:small}, one has that $\mathcal{I}_{U(1/\sqrt{2})} = \mathcal{S}_{U(1/\sqrt{2})}$, the two selection matrices shown in~\eqref{eq:selection matrix}. The following lemma shows that this relation holds in general, and thus the definition of $\mathcal{I}_U$ gives a simpler characterization of $\mathcal{S}_U$.
\begin{lemma} \label{lemma:suiu}
Consider an orthonormal matrix $U\in \R^{n \times r}$ with $r<n$. Then 
$\mathcal{S}_U=\mathcal I_U.$
\end{lemma}
\begin{proof}
 See Appendix~\ref{Appendix:Q-DEIM without tie-breaking rule}.
\end{proof}

\section{Theoretical foundations} \label{sec:theory}

The goal of this section is to provide theoretical foundations for combining DLRA with DEIM. As explained in the introduction, the index selection mechanism introduces discontinuities into the differential equations whenever indices are modified. Let us stress that this is not an exceptional event and bound to take place in the context of DLRA when, e.g., different regions of the spatial domain become relevant in the course of integrating over the time interval. 
We will construct a new continuous-time model based on differential inclusions to combine DEIM with a differential equation. 
Our framework and the discussion in this section applies to any DEIM procedure; the only exception is Section~\ref{sec:qdeim}, which works out an explicit expression for the differential inclusions when QDEIM is used.

\subsection{Oblique DEIM tangent space projection}

Given a rank-$r$ matrix $Y = U \Sigma V^\top \in \Mr \subset \Rmn$ with orthonormal $U \in \Rmr, V \in \Rnr$, we consider the oblique projections
\begin{equation} \label{eq: oblique projections range and corange}
    P_U^{\angle} = U (S_U^\top U)^{-1} S_U^\top, \qquad P_V^{\reflectangle} = S_V (V^\top S_V)^{-1} V^\top,
\end{equation}
where $S_U, S_V$ are selection matrices such that $S_U^\top U$ and $S_V^\top V$ are square and invertible.
These two projections define the linear matrix operator
\begin{equation} \label{eq: tangent oblique projection}
    \PPDEIM_Y\colon \R^{m\times n} \to T_Y \Mr, \quad Z \mapsto P_U^{\angle} Z - P_U^{\angle} Z P_V^{\reflectangle} + Z P_V^{\reflectangle}.
\end{equation}
One easily verifies the relation $\mathcal{P}_Y[\PPDEIM_Y Z]= \PPDEIM_Y Z$ for the orthogonal tangent space projection $\mathcal{P}_Y$ defined in~\eqref{eq: orthogonal tangent projection},
which implies that $\PPDEIM_Y Z\in \mathcal T_Y \Mr$ indeed holds. Similarly, elements on the tangent space are preserved, $\PPDEIM_Y X = X$ for every $X \in T_Y\Mr$. These two properties make 
$\PPDEIM_Y$ an oblique projection onto the tangent space $T_Y \Mr$.
In addition, $\PPDEIM_Y$ has an interpolatory property; it preserves matrix entries at the selected indices: $$S_U^\top \PPDEIM_Y[Z] S_V = S_U^\top Z S_V.$$

The following lemma bounds how much the error for the oblique tangent space projection $\PPDEIM_Y$ deviates from the best one, that is, the orthogonal projection. 
\begin{lemma}
\label{thm: deim_bound_2d}
    For $Y=U\Sigma V^\top\in \mathcal{M}_r$ with orthonormal $U \in \Rmr, V \in \Rnr$, suppose that $\|(S_U^\top U)^{-1}\|_2 \le C_{m,r}$ and 
    $\|(S_V^\top V)^{-1}\|_2 \le C_{n,r}$. Then the oblique projection $\PPDEIM_Y$ defined in~\eqref{eq: tangent oblique projection} satisfies
    $$\|Z - \mathcal{P}^{\angle}_Y Z\|_F \leq  C_{m,r} \cdot  C_{n,r} \cdot  \|Z -\mathcal{P}_Y Z\|_F, \quad \forall Z \in \R^{m\times n}.$$
\end{lemma}
\begin{proof}
Using their definitions, we can rewrite the projections as
\[
 Z - \mathcal{P}^{\angle}_Y Z = (I - P_U^{\angle}) Z (I - P_V^{\reflectangle}), \quad 
 Z - \mathcal{P}_Y Z = (I-UU^\top) Z (I-VV^\top).
\]
Combined with the identities
\begin{equation} \label{eq:projectors}
\begin{array}{lll}
 P_U^{\angle} UU^\top = UU^\top, \ & (I - P_U^{\angle}) UU^\top = 0, \ & (I - P_U^{\angle}) (I-UU^\top) = I - P_U^{\angle}, \\
 VV^\top P_V^{\angle}  = VV^\top, \ &
 VV^\top (I-P_V^{\reflectangle}) = 0,\ & 
 (I-VV^\top)(I-P_V^{\reflectangle}) = I-P_V^{\reflectangle},
\end{array}
\end{equation}
this implies 
\[
    \|Z - \mathcal{P}^{\angle}_Y Z\|_F = \|  (I - P_U^{\angle})  ( Z - \mathcal{P}_Y Z )(I-P_V^{\reflectangle})\|_F 
    \le \| I - P_U^{\angle} \|_2 \| I-P_V^{\reflectangle}\|_2  \|Z - \mathcal{P}_Y Z \|_F.
\]
The proof is concluded by noting that
\[
 \| I - P_U^{\angle} \|_2 = \| P_U^{\angle} \|_2  =\|(S_U^\top U)^{-1}\|_2,\quad  \| I-P_V^{\reflectangle}\|_2   = \| P_V^{\reflectangle}\|_2  = \|(S_V^\top V)^{-1}\|_2.
\]
\end{proof}

For later purposes, we note that the result of Lemma~\ref{thm: deim_bound_2d} also implies
\begin{equation} \label{eq:boundproj}
 \|\PPDEIM_Y\|_2 = \|\mathcal{I} - \PPDEIM_Y\|_2 \le C_{m,r} \cdot  C_{n,r},´
\end{equation}
where $\|\PPDEIM_Y\|_2$ denotes the norm of the linear operator $\PPDEIM_Y$ induced by the Frobenius norm on $\R^{m\times n}$ and $\mathcal{I}$ is the identity on $\R^{m\times n}$.

Table~\ref{tab:algo_comparison} gives an overview of the upper bounds on $C_{m,r}$ and $C_{n,r}$ (when replacing $m$ by $n$) for the various DEIM procedures. For example, when
Osinsky’s method~\cite{osinsky_close_2023} or ARP~\cite{cortinovis_adaptive_2024} are used then the bound of Lemma~\ref{thm: deim_bound_2d} holds with
$C_{m,r} \cdot  C_{n,r} \le \sqrt{1+r(m-r)} \sqrt{1+r(n-r)}$ (in expectation).

\subsection{Definition and properties of DLRA-DEIM}  \label{sec: DI}

We now aim at deriving a suitable continuous-time model that incorporates the oblique projection $\PPDEIM_Y$ into DLRA.
As already emphasized, the naive way of simply replacing orthogonal by oblique projections in DLRA leads to a differential equations~\eqref{eq:dlardeimintro} that is ill-posed in the presence of the discontinuities caused by DEIM. In the literature on differential equations, such discontinuities are commonly resolved by considering Filippov inclusions; see \cite{ acary2008numerical,Smirnov2002Introduction} for an introduction.
In particular, we will follow the framework described in~\cite[Section 4.4]{Smirnov2002Introduction}.

Let us first introduce 
$$f(Y):=\mathcal{P}^\angle_{Y}[F(Y)],$$
which we assume to be always well-defined (for example, because QDEIM is used with tie-breaking). Letting $P(\Rmn)$ denote the power set of $m\times n$ matrices, we then define the set-valued map $\mathcal{F}: \mathcal{M}_r\rightarrow {P}( \mathbb{R}^{m\times n})$ by
\begin{equation} \label{eq: filippov inclusion}
    \mathcal{F}(Y):=\bigcap_{\epsilon>0}\overline{\text{conv}}\big[f\big((Y+\epsilon \mathcal{B})\cap \mathcal{M}_r\big)\big],
\end{equation}
 where $\mathcal{B} = \{X\in \mathbb{R}^{m\times n}:\|X\|_F\leq 1\}$ denotes the ball of radius $1$ and we use set-valued notation:
 $$f\big((Y+\epsilon\mathcal{B})\cap \mathcal{M}_r\big) = \bigcup_{Z\in (Y+\epsilon\mathcal{B})\cap \mathcal{M}_r}f(Z).$$ 
The notation $\text{conv}(\cdot)$ denotes the convex hull of a set and $\overline{\text{conv}}(\cdot)$ denotes its closure. 
The set $\mathcal{F}(Y)$ is always non-empty because $f(Y)\in \mathcal{F}(Y)$. When $f$ is continuous at $Y$, then $f(Y)$ is the only element in $\mathcal{F}(Y)$; see, e.g.,~\cite[P. 105]{Smirnov2002Introduction}.

Consider an absolute continuous function $Y:[0,T]\rightarrow \mathcal{M}_r$, which is differentiable almost everywhere on $[0,T]$. Generalizing~\eqref{eq: DLRA}, we say that \emph{$Y(t)$ is a solution to DLRA-DEIM} on an interval $[0,T]$ if it satisfies the \emph{differential inclusion}
\begin{equation}
    \label{eq:DI}
    Y(0)=Y_0\in \mathcal{M}_r\quad \text{and}\quad \dt{Y}(t)\in \mathcal{F}(Y(t)) \quad \text{for almost every} \quad t\in [0,T].
\end{equation}
Let us recall that $\dt{Y}(t)=f(Y(t))$ whenever $f$ is continuous at $Y(t)$, in which case the inclusion coincides with~\eqref{eq:dlardeimintro}.

To provide more insight into the abstract expression~\eqref{eq:DI}, we first provide a lemma that allows us to reformulate $\mathcal F(Y)$ and shows that all values returned by $\mathcal{F}$ are contained in $T_Y\Mr$.
\begin{lemma} \label{lemma: limit_f}
Suppose that Assumption~\ref{ass: bounded F} holds. Then, for $Y\in \M_r$ it holds that
\begin{equation}
\label{eq:limit_f}
     \mathcal{F}(Y)=\textup{conv}\{X\in \mathbb{R}^{m\times n} \mid \exists\{Y_i\}\subseteq \M_r \text{ s.t. } Y_i\rightarrow Y \text{ and } f(Y_i)\rightarrow X\}\subseteq T_Y\Mr.
\end{equation}   
\end{lemma}
\begin{proof}
The equality in~\eqref{eq:limit_f} is standard in the literature on differential inclusions, see, e.g.,~\cite[Lemma 1]{Paden1987calculus}. For the convenience of the reader, details on the proof specific to our setup are provided in Appendix~\ref{App: supplementary lemma}.

It remains to prove the inclusion $\mathcal{F}(Y)\subseteq T_Y\Mr$.
Choose some $X\in \R^{m\times n}$ such that there exists $\{Y_i\}\subseteq \M_r$ with $Y_i\rightarrow Y$ and $f(Y_i)\rightarrow X$. Because the ranks of $Y_i,Y$ remain constant, one can choose continuous singular value decompositions $Y_i=U_i\Sigma_iV_i^\top$ and $Y=U\Sigma V^\top$ such that
$U_i\to U$ and $V_i \to V$.  Because of $f(Y_i)\in T_{Y_i}\mathcal{M}_r$, we obtain that
\begin{align*}
0 &= f(Y_i)-\mathcal{P}_{Y_i}[f(Y_i)]=f(Y_i)-(U_iU_i^\top f(Y_i)+f(Y_i)V_iV_i^\top-U_iU_i^\top f(Y_i)V_iV_i^\top) \\
&\to f(Y)-(U U^\top f(Y)+f(Y)V V^\top-U U^\top f(Y)VV^\top) = f(Y)-\mathcal{P}_{Y}[f(Y)].
\end{align*} 
This implies $X\in T_{Y}\mathcal{M}_r.$ The proof is completed by noting that the convex hull of all such $X$ remains in 
$T_{Y}\mathcal{M}_r$. 
\end{proof}

\subsubsection{Structure of $\mathcal{F}$ for QDEIM} \label{sec:qdeim}

Lemma~\ref{lemma: limit_f} remains an abstract characterization of $\mathcal{F}(Y)$. 
To obtain a concrete expression, we now consider the case when QDEIM (Algorithm \ref{alg:QDEIM}) is used for index selection. Let us recall that the set $\mathcal I_U$ defined in~\eqref{eq:allselections} collects all feasible selection matrices returned by QDEIM without tie-breaking. The following lemma shows that 
$\mathcal{F}(Y)$ is the convex polytope spanned by the (finitely many different) oblique projections of $F(Y)$ that arise from any such feasible choice of selection matrices. 
\begin{lemma} \label{lemma:qdeimf}
    Consider $Y=U\Sigma V^\top\in \Mr$ with orthonormal $U \in \Rmr$, $V \in \Rnr$ for $r < m,n$. Suppose that $F$ is continuous and that Assumption~\ref{ass: bounded F} holds. If QDEIM
    is used to define the oblique projection $\mathcal P_Y^\angle$ in~\eqref{eq: tangent oblique projection} then
    \begin{equation*}
        \mathcal{F}(Y)= \textup{conv}\big\{P_U^{\angle} F(Y) - P_U^{\angle} F(Y) P_V^{\reflectangle} + F(Y) P_V^{\reflectangle}\big|S_U\in \mathcal{I}_U\text{ and } S_V\in \mathcal{I}_V \big\},
    \end{equation*}
    with the oblique projections $P_U^{\angle}$, 
    $P_V^{\reflectangle}$ defined as in~\eqref{eq: oblique projections range and corange} depending on $S_U, S_V$.
\end{lemma}
\begin{proof}
We will make use of the fact that the sets $\mathcal I_U, \mathcal I_V$ as well as the elements picked by QDEIM with tie breaking are invariant under the choice of orthonormal bases $U,V$; see~\cite[Theorem 2.1]{drmac2016new}.

    To prove the inclusion ``$\subseteq$", let $X\in \R^{m\times n}$ be such that
    there exists a sequence $\{Y_i\}\subseteq \M_r$ with $Y_i\rightarrow Y$ and $f(Y_i)\rightarrow X$. As in the proof of Lemma~\ref{lemma: limit_f}, we can choose singular value decompositions $Y_i=U_i\Sigma_iV_i^\top$ and $Y=U\Sigma V^\top$ such that
$U_i\to U$ and $V_i \to V$. Because the set of selection matrices is a finite set, we can take subsequences to ensure $S^{\mathsf{Q}}_{U_i} \to S_U$, 
$S^{\mathsf{Q}}_{V_i} \to S_V$ for some selection matrices $S_U, S_V$. By the definition~\eqref{eq:vicinity}, we have 
that $S_U\in \mathcal S_U$, $S_V\in \mathcal S_V$.
Lemma~\ref{lemma:suiu} implies $S_U\in \mathcal I_U$, $S_V\in \mathcal I_V$. Because these are feasible choices of QDEIM without tie breaking, it follows that 
$P_{U}^{\angle}, P_{V}^{\reflectangle}$ are well-defined and 
$P_{U_i}^{\angle}\to P_{U}^{\angle}$,
$P_{V_i}^{\reflectangle} \to P_{V}^{\reflectangle}$. Together with the continuity of $F$, we therefore obtain that
     $$X=\lim_{i\rightarrow \infty} f(Y_i)=\lim_{i\rightarrow \infty} \big(P_{U_i}^{\angle} F(Y_i) - P_{U_i}^{\angle} F(Y_i) P_{V_i}^{\reflectangle} + F(Y_i) P_{V_i}^{\reflectangle}\big)=P_{U}^{\angle} F(Y) - P_{U}^{\angle} F(Y) P_V^{\reflectangle} + F(Y) P_V^{\reflectangle}.$$
Taking convex hulls and using the characterization of $\mathcal F(Y)$ by
    Lemma~\ref{lemma: limit_f} completes the proof of the inclusion ``$\subseteq$''.

For the other inclusion ``$\supseteq$", let $S_U \in \mathcal I_U$ and $S_V \in \mathcal{I}_V$. By Lemma~\ref{lemma:suiu} and the definition~\eqref{eq:vicinity}
of $\mathcal S_U, \mathcal S_V$, there exist orthonormal matrix sequences 
$U_i\rightarrow U$, $V_i\rightarrow V$ such that $S^\mathsf{Q}_{U_i}\rightarrow S_U$, $S^\mathsf{Q}_{V_i}\rightarrow S_V$.
Setting $Y_i:=U_i\Sigma V^\top_i\in \Mr$ gives $Y_i\rightarrow Y=U\Sigma V^\top$. Using the same reasoning as above, the limit $\lim_{i\rightarrow \infty} f(Y_i)$ exists and is in $\mathcal{F}(Y)$. Again, taking convex hulls and using 
Lemma~\ref{lemma: limit_f} completes the proof of ``$\supseteq$".
\end{proof}

\begin{figure}[ht!]
\small
\centering
\includegraphics[scale=0.4]{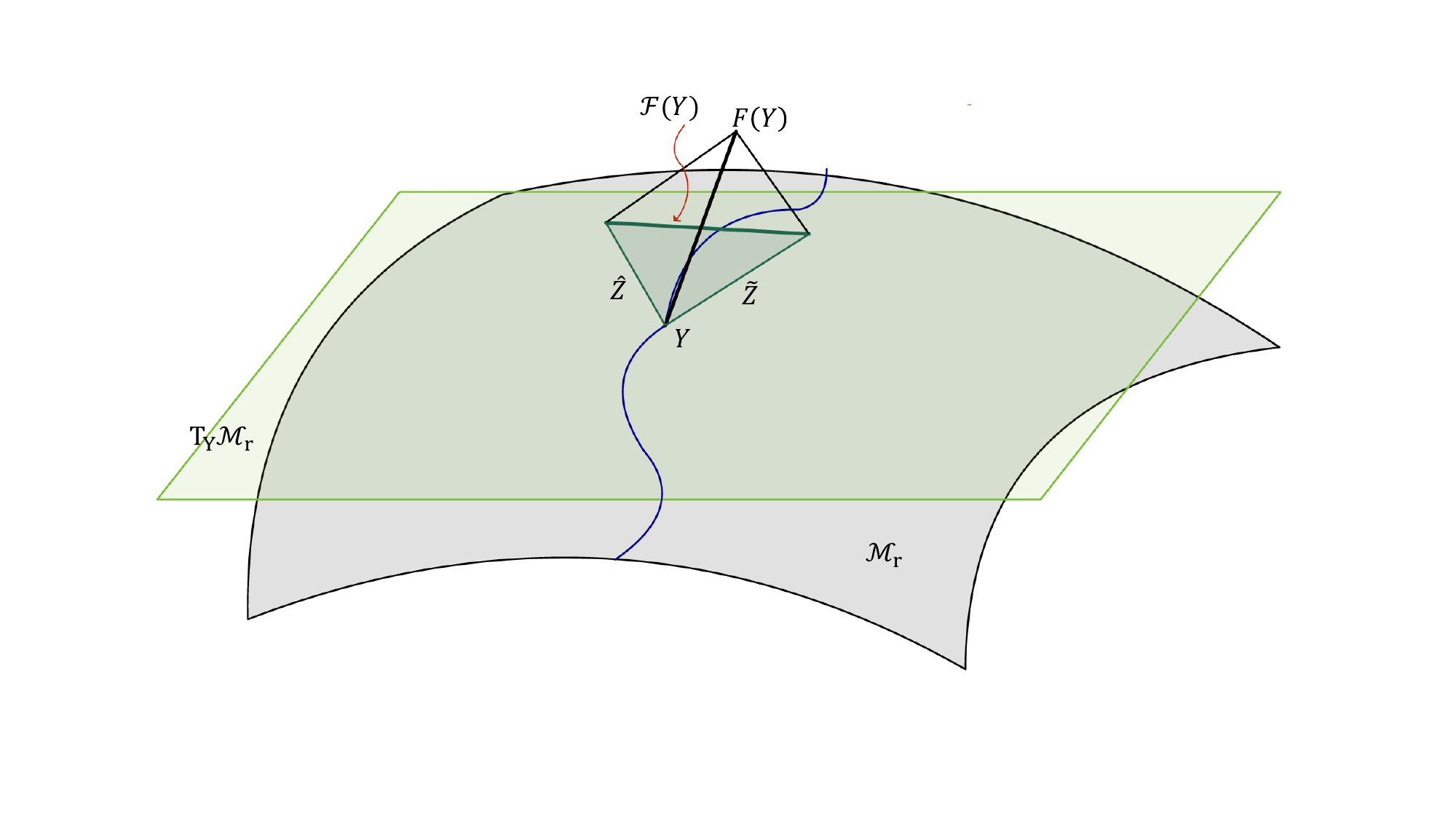}
\caption{Illustration of the set $\mathcal F(Y)$ in DLRA-DEIM~\eqref{eq:DI} 
with QDEIM for the matrix $Y$ in~\eqref{eq:smallY}.}
\label{fig: DLRA-DEIM}
\end{figure}

Let us illustrate the set $\mathcal{F}(Y)$ for the following rank-$2$ matrix: 
\begin{equation} \label{eq:smallY}
Y= \underbrace{\begin{pmatrix}1/\sqrt{2} &  0 \\ 1/\sqrt{2}& 0 \\ 0 & 1 \end{pmatrix}}_{:=U}\underbrace{\begin{pmatrix} 2 &  0 \\ 0& 1\end{pmatrix}}_{:=\Sigma}\underbrace{\begin{pmatrix} 1/2 &  \sqrt{3/4} &0\\ 0&0& 1\end{pmatrix}}_{:=V^\top}. 
\end{equation}
As explained in Example~\ref{ex:small}, the set $\mathcal I_U$ consists of the the two selection matrices 
$$\tilde{S}_U = \begin{pmatrix} 0 & 0 \\ 0 & 1 \\ 1 & 0 \end{pmatrix}, \quad \hat{S}_U = \begin{pmatrix} 0 & 1 \\ 0 & 0 \\ 1 & 0 \end{pmatrix}.$$
It can be easily verified that the set $\mathcal I_V$ only contains a single selection matrix $S_V$. By Lemma~\ref{lemma:qdeimf},
$$\mathcal{F}(Y)=\text{conv}\{ \tilde{Z},\hat{Z}\}=\{\alpha \tilde{Z}+(1-\alpha)\hat{Z}|\alpha\in [0,1]\} \subset T_Y \Mr,$$
with \begin{align*}\tilde{Z} &= U(\tilde{S}_U^\top U)^{-1}\tilde{S}_U^\top F(Y) - U(\tilde{S}_U^\top U)^{-1}\tilde{S}_U^\top F(Y) P_V^{\reflectangle} + F(Y) P_V^{\reflectangle} \\ \hat{Z} &=  U(\hat{S}_U^\top U)^{-1}\hat{S}_U^\top F(Y) - U(\hat{S}_U^\top U)^{-1}\hat{S}_U^\top F(Y) P_V^{\reflectangle}+ F(Y) P_V^{\reflectangle}.\end{align*}
In other words, $\mathcal{F}(Y)$ is a straight line in the tangent space as illustrated in Figure~\ref{fig: DLRA-DEIM}.

\subsection{Analysis of DLRA-DEIM}

In this section, we will analyze DLRA-DEIM~\eqref{eq:DI}, establishing an exactness result as well as existence and approximation properties of the solution.
We first show that DLRA-DEIM inherits an exactness property of DLRA, that is, it recovers solutions of the~\eqref{eq: full order model}. 
\begin{lemma}[Exactness property]
Suppose that the solution of the full-order model~\eqref{eq: full order model} 
is unique and satisfies 
$A(t) \in \Mr$ for every $t \in [0,T]$. Then $A(t)$ is a solution to the differential inclusion~\eqref{eq:DI}. If $F$ is continuous, then this is the unique differentiable solution to~\eqref{eq:DI}, provided that QDEIM is used for index selection.
\end{lemma}
\begin{proof}
    The assumptions imply that $A(0) \in \Mr$ and $\dt{A}(t) \in T_{A(t)} \Mr$ for every $t \in [0,T]$. By Lemma~\ref{thm: deim_bound_2d}, the oblique projection $\PPDEIM_{A(t)}$ preserves the tangent space:
    $$\dt{A}(t) = F(A(t)) = \PPDEIM_{A(t)} F(A(t)) \in \mathcal F(A(t)),$$
    which shows that $A(t)$ verifies the desired differential inclusion. If QDEIM is used then Lemma~\ref{lemma:qdeimf} shows that every element of $\mathcal F(A(t))$ is given by such an oblique projection (for possibly other selection matrices $S_U, S_V$), which preserves the tangent space and thus equals $F(A(t))$. In other words, the differential inclusion coincides with the full-order model, which shows that its unique differentiable solution is given by $A(t)$.
\end{proof}

When $A(t)$ is not in $\Mr$, the existence of a solution to the differential inclusion~\eqref{eq:DI} is less obvious. 
The next theorem establishes the existence on an interval $[0, t_{\max})$, where the lower bound on $t_{\max}$ depends on the boundedness of $F$, the distance of the initial value $Y(0)$ from the boundary of $\mathcal M_r$ (measured in terms of its smallest nonzero singular value $\sigma_{r}(Y(0))$), as well as the constants $C_{m,r}$, $C_{n,r}$ governing the quality of the DEIM procedure.

\begin{theorem} \label{thm: DI_existent}
 Suppose that Assumption~\ref{ass: bounded F} holds. Then there exists a solution $Y(t)\in \mathcal{M}_r$ to the differential inclusion~\eqref{eq:DI} for $t\in [0,t_{max}),$ with $t_{\max} > t_1: =\frac{\sigma_{r}(Y(0))}{2 C_{m,r} C_{n,r} B}$, where  
 $C_{m,r}$, $C_{n,r}$ are as in Lemma~\ref{lemma: deviation of P_diem}. 
\end{theorem}
\begin{proof}
The proof follows by applying a general result on differential inclusions~\cite[Corollary 1.1]{Haddad1981Monotone}, which is stated for $\M_r$ as Theorem~\ref{thm: existant} in Appendix~\ref{App: DI}. To apply the theorem we need to first verify its assumptions on $\mathcal{F}$.

It immediately follows from Lemma~\ref{lemma: F_usc} and the definition of $\mathcal{F}$ 
that $\mathcal{F}$ is upper semicontinuous and that it returns non-empty, convex and closed sets. Moreover, for any $X\in f((Y+\epsilon \mathcal{B})\cap \M_r)$ 
it holds that
    \begin{equation}\label{eq:F bound}
    \|X\|_F\leq \|\PPDEIM_Y\|_2 \cdot B \leq C_{m,r}C_{n,r} B,\end{equation}
    where we used Assumption~\ref{ass: bounded F} and~\eqref{eq:boundproj}. Therefore, $f((Y+\epsilon \mathcal{B})\cap \M_r)$ is always a bounded set. As intersection and convex closure preserve boundedness, we obtain that $\mathcal{F}(Y)$ is bounded by the right-hand side of~\eqref{eq:F bound}. Together with their closedness, this implies that the values of $\mathcal{F}$ are compact. 

Because $f(Y)$ is contained in $\mathcal{F}(Y)\cap T_Y \mathcal{M}_r$, this set is non-empty for every $Y\in \mathcal{M}_r$. Thus, all assumptions of Theorem~\ref{thm: existant} hold, which implies the existence of a solution $Y(t)\in \M_r$ for $t\in[0,t_{\max})$. To state a lower bound for $t_{\max}$, we first note that 
$\sigma_r(Y(0))>0$ because $Y(0)\in \mathcal{M}_r$ and, thus, the set $$\{Y\in \mathcal{M}_r \mid \|Y-Y_0\|_F\leq \sigma_r(Y(0))/2\}\subseteq \mathcal{M}_r$$ is closed. By~\eqref{eq:F bound}, $\mathcal{F}(Y)$ is bounded by $C_{m,r}C_{n,r}B.$ Hence, Theorem~\ref{thm: existant} implies the lower bound $t_{\max}>\frac{\sigma_{r}(Y(0))}{2 C_{m,r}C_{n,r}B}$.
\end{proof}

The construction in the proof of Theorem \ref{thm: DI_existent} ensures that $\sigma_{r}(Y(t_1)) > 0$ and, hence, we can apply the theorem again to extend the solution to a larger time interval $[0, t_2]$ with $t_2 > t_1$. This argument can be repeated but the maximal attainable interval $[0, T)$ remains bounded
when the solution reaches the boundary of the manifold $\mathcal M_r$, at which point $\sigma_{r}(Y(t))$ converges to zero. This limitation, caused by rank deficiencies, is common to continuous-time DLRA models.

We now aim at quantifying the error between a solution to the DLRA-DEIM and the solution to the full-order model. We start by extending the property of Assumption~\ref{ass: epsilon close}.

\begin{corollary} \label{lemma: deviation of P_diem}
Let $F$ be continuous.    Suppose that Assumptions~\ref{ass: epsilon close} and~\ref{ass: bounded F} hold. Then for any $Y \in \Mr$ in a suitable neighborhood of the solution to~\eqref{eq: full order model},
    we have that
$$\|\dot{Y}-F(Y)\|_F\leq C_{m,r} C_{n,r} \cdot \varepsilon_r \quad \forall \dot{Y}\in \mathcal{F}(Y).$$
\end{corollary}
\begin{proof}
 By Lemma~\ref{lemma: limit_f}, 
 \begin{equation*}
     \dot{Y}\in \mathcal{F}(Y)= \textup{conv}\{X\in \mathbb{R}^{m\times n} \mid \exists\{Y_i\}\subseteq \M_r \text{ s.t. } Y_i\rightarrow Y \text{ and } f(Y_i)\rightarrow X\}.
 \end{equation*} Take any such $X$, for which there exists a sequence $\{Y_i\}\subseteq\M_r$ with $Y_i\rightarrow Y$ 
 and $f(Y_i)\rightarrow X$. By the continuity of $F$, 
 \[
     \|X-F(Y)\|_F = \|\lim_{i\rightarrow \infty}(f(Y_i)-F(Y_i))\|_F = 
     \lim_{i\rightarrow \infty}\| \PPDEIM_{Y_i}[F(Y_i)] -  F(Y_i)\|_F
\leq C_{m,r} C_{n,r} \cdot \varepsilon_r,
 \] where we have used Assumption~\ref{ass: epsilon close} and Lemma~\ref{thm: deim_bound_2d}. The statement follows by writing $\dot{Y} - F(Y)$ as a convex combination of such $X-F(X)$
 and applying the bound to every summand in the combination.
 \end{proof}
Finally, we provide an error bound for DLRA-DEIM, in a form similar to Theorem~\ref{thm: DLRA_acc}.
\begin{theorem} \label{thm: dlra-deim error}
Suppose that Assumptions \ref{ass: one-sided F}, \ref{ass: epsilon close},and \ref{ass: bounded F} hold, and let
$A(t)$ be the solution of the differential equation 
$$\dt{A}(t)=F(A(t)),\quad A(0)=A_0\in \mathbb{R}^{n \times n}.$$
Then for any solution $Y(t) \in \M_r$ of \eqref{eq:DI} defined on an interval $[0,T]$, we have the error bound
     $$\|A(t)-Y(t)\|_F\leq e^{\ell t}\|A(0)-Y(0)\|_F+C_{m,r} C_{n,r} \cdot \varepsilon_r\int^{t}_0e^{\ell(t-s)} \mathrm{d}s.$$
 \end{theorem}
  \begin{proof}
  This result is proven by following the proof of the corresponding DLRA result~\cite[Lemma 1]{kieri2019projection} verbatim with $\varepsilon_r$ replaced by $C_{m,r} C_{n,r}\cdot \varepsilon_r$.
 \end{proof}

\section{PRK-DEIM: Projected Runge--Kutta DEIM} \label{sec: PRK-DEIM}

In this section, we introduce a family of integrators for DLRA-DEIM, called projected Runge--Kutta DEIM (PRK-DEIM) methods. They are a generalization of the PRK methods from \cite{kieri2019projection} for DLRA, using oblique instead of orthogonal tangent space projections.

Both, PRK and PRK-DEIM are built on classical explicit Runge--Kutta (RK) methods~\cite{hairersolving,hairer2006numerical,hairer2003geometric}. Applied to the full-order model~\eqref{eq: full order model} with step size $h>0$, 
an \textit{explicit RK method with $s$ stages} produces the approximations $A_i \approx A(ih)$ by successively computing
\begin{equation} \label{eq: Runge--Kutta method}
\begin{aligned}
\tilde{Z}_j &= A_i + h \sum_{\ell = 1}^{j-1} a_{j \ell} F(\tilde{Z}_{\ell}), \quad j=1, \ldots, s, \\
A_{i+1} &= A_i + h \sum_{j=1}^s b_j F(\tilde{Z}_j).
\end{aligned}
\end{equation}
Here, the coefficients $a_{ij}$ and $b_j$ for $i,j=1,\ldots,s$ determine the specific choice of RK method. Such explicit methods are well suited for non-stiff problems and the flexibility of having several stages can be used to design high-order methods, that is, the approximation error is $O(h^p)$ for some $p\ge 2$, provided that $F$ is sufficiently smooth.

Applying classical RK methods directly to the DLRA~\eqref{eq: DLRA} comes with the major disadvantage that the iterates $A_i$ in general do \emph{not} remain on the manifold $\mathcal M_r$, which annihilates all computational advantages gained from working with low-rank matrices. Even when replacing $A_i$ by an approximation on $\mathcal M_r$ through rank-$r$ truncation, the intermediate quantities appearing in~\eqref{eq: Runge--Kutta method} might have large or even full rank. The PRK methods from~\cite{kieri2019projection} partly address this problem by additionally applying orthogonal tangent space projections to the stages, limiting the ranks of intermediate RK stages to $2r$. These methods are proven to be robust to the presence of small singular values and preserve the original RK convergence rate up to order three. However, for general nonlinearities PRK methods may still require to repeatedly evaluate the complete $m\times n$ matrix $F(\cdot)$ explicitly, which can become a major computational bottleneck. 
Motivated by DLRA-DEIM~\eqref{eq:DI}, we propose to address this limitation by using oblique tangent space projection.

For the Euler method, the corresponding projected Euler-DEIM method takes the form
\begin{equation} \label{eq: projected Euler DEIM}
    Y_0 = \Tr(A_0), \quad Y_{i+1}=\Tr \big(Y_i+h\mathcal{P}^\angle_{Y_i}[F(Y_i)]\big), \quad i = 0,1,\ldots.
\end{equation}
Here, $\Tr : \Rmn \rightarrow \Mr$ denotes rank-$r$ truncation, which returns a best rank-$r$ approximation of a matrix in the Frobenius norm.
We represent each $Y_i\in \mathcal M_r$ in factorized form $Y_i = U_i S_i V_i^\top$ for orthonormal matrices $U_i \in \Rmr$, $V_i \in \Rnr$.
Recalling that
$$\mathcal{P}^\angle_{Y_i}[F(Y_i)] = P_{U_i}^{\angle} F(Y_i) - P_{U_i}^{\angle} F(Y_i) P_{V_i}^{\reflectangle} + F(Y_i) P_{V_i}^{\reflectangle},$$ 
and the definition~\eqref{eq: oblique projections range and corange} of the DEIM oblique projections $P_{U_i}^{\angle}$ and $P_{V_i}^{\reflectangle}$, we see that
the computation of this projection (in factorized form) only requires the evaluation of $r$ rows and $r$ columns of~$F$, corresponding to the row indices selected by DEIM for $U_i$ and $V_i$.
This is usually much cheaper than evaluating the whole matrix $F(Y_i)$. For example, when $F$ is defined element-wise (e.g., a cubic nonlinearity) only $O((m+n)r^2)$ operations are required.
Note that the truncation operator $\Tr$ in~\eqref{eq: projected Euler DEIM} is applied to a matrix of rank at most $2r$, which requires another $O((m+n)r^2)$ operations.

The extension from~\eqref{eq: projected Euler DEIM} to higher-order RK methods is straightforward. Given a  classical RK method~\eqref{eq: Runge--Kutta method}, we define the corresponding \textit{PRK-DEIM method} as
\begin{equation} \label{eq: PRK-DEIM}
\begin{aligned}
Z_j &= Y_i + h \sum_{\ell = 1}^{j-1} a_{j \ell} \mathcal P_{\Tr(Z_{\ell})}^{\angle}\left[F(\Tr(Z_{\ell}))\right], \quad j=1, \ldots, s, \\
Y_{i+1} &= \Tr \Big(Y_i + h \sum_{j=1}^s b_j \mathcal P_{\Tr(Z_j)}^{\angle}\left[F(\Tr(Z_j))\right] \Big).
\end{aligned}
\end{equation}
Similar to PRK methods, the rank growth in~\eqref{eq: PRK-DEIM} is limited to at most $2sr$ and the rank-$r$ truncation with $\Tr$ can be performed efficiently. 
As explained above for the projected Euler-DEIM method, the task of evaluating the projected terms $\mathcal P_{\Tr(Z_i)}^{\angle}\left[F(\Tr(Z_i))\right]$, $i = 1, \ldots, s$, reduces to the evaluation of $r$ selected columns and rows of $F(\Tr(Z_i))$ have only to be evaluated on the interpolation points. To select these columns and rows, the DEIM procedure needs to be run $2s$ times in every time step. It is therefore essential to limit the cost of such a selection procedure.

\subsection{Error analysis} \label{sec:error}

We now aim at establishing approximation results for~\eqref{eq: PRK-DEIM}. There are two main error sources: (1) restricting the (discrete) dynamics to the tangent space of $\mathcal M_r$, and (2) performing time discretization. The first error is largely due to the modeling error $\varepsilon_r$ of DLRA, which is captured by Assumption~\ref{ass: epsilon close}. One therefore expects this error to be of order $\varepsilon_r$ (see Theorem~\ref{thm: DLRA_acc}), multiplied by $C_{m,r}C_{n,r}$ to take into account the use of oblique projections (see Lemma~\ref{thm: deim_bound_2d}).
Following~\cite{kieri2019projection}, we do not study the second error separately (that is, the convergence of the time discretization~\eqref{eq: PRK-DEIM} to DEIM-DLRA~\eqref{eq:DI}), but directly compare~\eqref{eq: PRK-DEIM} with the full-order model. Since the proofs follow very closely those for the PRK methods in~\cite{kieri2019projection}, we will be brief and only point out the essential differences.

We first present results and proofs for the 
projected Euler-DEIM method~\eqref{eq: projected Euler DEIM}; this allows us to highlight the differences to~\cite{kieri2019projection} without clutter.

\begin{lemma}[Local error of projected Euler-DEIM] \label{lemma: local error projected Euler DEIM} Suppose that Assumptions \ref{ass: Lipschitz}, \ref{ass: one-sided F},~\ref{ass: epsilon close} hold. Let $\phi_F^t$ denote the exact flow for the full-order model, that is, $\phi_F^t(A)$ is the solution of~\eqref{eq: full order model} at time $t$, started at~$A~\in~\Rmn$.
Then the local error of the projected Euler-DEIM method \eqref{eq: projected Euler DEIM} is bounded by
    $$\big\|Y_{i+1} - \phi_F^h(Y_i) \big\|_F \leq C h (h + C_{m,r}C_{n,r} \cdot \varepsilon_r),$$
    where $C$ is the same constant as in~\cite[Lemma 3]{kieri2019projection}. In particular, $C$ depends on the Lipschitz constant $L$ of $F$ but it does not depend on $h$ or $\varepsilon_r$.
\end{lemma}

\begin{proof}
    The proof is along the lines of the proof of~\cite[Lemma 3]{kieri2019projection}. The main difference is that the orthogonal projection $\mathcal P_{Y_i}$ is replaced by the oblique projection $\mathcal  P_{Y_i}^{\angle}$ in the identity
    $$Y_i + h \mathcal P_{Y_i}^{\angle} F(Y_i) = Y_i + h F(Y_i) + h (\mathcal{I}-\mathcal P_{Y_i}^{\angle})F(Y_i).$$
    Lemma~\ref{thm: deim_bound_2d} gives a bound for the third term: $\norm{h (\mathcal{I}-\mathcal P_{Y_i}^{\angle})F(Y_i)}_F \leq C_{m,r} C_{n,r} \cdot h \cdot \varepsilon_r$. The first two terms constitute the classical Euler method, which satisfies the local error bound:
    $$\big\| \phi_F^h(Y_i) - Y_i - h F(Y_i) \big\|_F \leq C_L h^2,$$
    where $C_L$ is a constant that depends on  $L$ but not on $h$.
    We can therefore regard projected-Euler DEIM as the truncation of a perturbed flow:
    $$Y_{i+1} = \Tr(\phi_F^h(Y_i) + h \Delta), \quad \text{ with } \norm{\Delta}_F = O( h + C_{m,r} C_{n,r} \cdot \varepsilon_r).$$
    The rest of the proof of \cite[Lemma 3]{kieri2019projection} holds verbatim with $\varepsilon = C_{m,r} C_{n,r} \, \varepsilon_r$.
\end{proof}

In the same way the global error bound for PRK~\cite[Theorem 4]{kieri2019projection} follows from~\cite[Lemma~3]{kieri2019projection}, the following result follows from
Lemma~\ref{lemma: local error projected Euler DEIM}. 

\begin{theorem}[Global error of projected Euler DEIM] \label{thm: global error projected Euler DEIM}
Suppose that Assumptions \ref{ass: Lipschitz}, \ref{ass: one-sided F},~\ref{ass: epsilon close} hold
and set $\delta := \|A_0 - Y_0 \|_F$.
Then there is $h_0 > 0$ such the projected Euler-DEIM method~\eqref{eq: projected Euler DEIM}
satisfies
$$\| A(N h) - Y_N \|_F \leq C(\delta + C_{m,r} C_{n,r} \cdot \varepsilon_r + h),$$
for every $0 < h < h_0$ and every $N \in \mathbb N$ such that $Nh \leq T$. Here, $C$ is the constant from~\cite[Theorem 4]{kieri2019projection}, which does not depend on $h$ or $\varepsilon_r$. 
\end{theorem}
\begin{proof}
    The proof of \cite[Theorem 4]{kieri2019projection} holds verbatim using the bound from Lemma \ref{lemma: local error projected Euler DEIM} and $\varepsilon = C_{m,r} C_{n,r}\, \varepsilon_r$.
\end{proof}

The analysis of higher-order methods requires the following extension of~\cite[Lemma 7]{kieri2019projection}; see Appendix \ref{App: proof of projected stages error} for the proof.

\begin{lemma}[Local error of PRK-DEIM stages] \label{lemma: projected stages error}
Suppose that Assumptions \ref{ass: Lipschitz}, \ref{ass: one-sided F}, \ref{ass: epsilon close} hold.
Consider the stages $\tilde Z_1, \ldots, \tilde Z_s$ of an explicit RK method~\eqref{eq: Runge--Kutta method}
and the stages $Z_1, \ldots, Z_s$ of the corresponding PRK-DEIM method~\eqref{eq: PRK-DEIM} with $A_i = Y_i \in \mathcal M_r$.
Let $q_1 \leq q_2 \leq \ldots \leq q_s$ denote the order of each stage of~\eqref{eq: Runge--Kutta method} and define $\tilde{q}_j = \min(q_j, q_2+1)$. Then the following bounds hold:
\begin{equation*}
\| Z_j - \tilde{Z}_j \|_F \leq \left\{\begin{aligned} &0 &&\text{ for } j=1, \\ &c_1 h \left( c_3 \cdot C_{m,r} C_{n,r} \cdot \varepsilon_r + h^{q_2 +1} \right) && \text{ for }
j=2, \ldots, s,
\end{aligned} \right.
\end{equation*}
\begin{equation*}
\big\| {\mathcal P_{\Tr(Z_j)}^{\angle} F(\Tr(Z_j)) - F(\tilde Z_j)} \big\|_F \leq \left\{ \begin{aligned} &C_{m,r} C_{n,r} \cdot \varepsilon_r && \text{ for } j=1, \\
&c_2 \left( c_4 \cdot C_{m,r} C_{n,r} \cdot \varepsilon_r + h^{\tilde q_j +1} \right) && \text{ for } j=2, \ldots, s. \end{aligned} \right.
\end{equation*}
The constants $c_1,\ldots, c_4$ do not depend on $h$ or $\varepsilon_r$.  
\end{lemma}

\begin{theorem}[Global error of PRK-DEIM] \label{thm: convergence of PRK-DEIM}
Considering the setting and assumptions of Lemma~\ref{lemma: projected stages error},
set $\delta := \|A_0 - Y_0 \|_F$ and
\begin{equation*}
q = \left\{
\begin{aligned}
&\min \{p, q_2 + 1\}, && \text{ if } b_2 \neq 0, \\
&\min \{p, q_3 + 1, q_2 +1 \}, && \text{ if } b_2 = 0.
\end{aligned} \right.
\end{equation*}
Then there is $h_0 > 0$ such that the PRK-DEIM method~\eqref{eq: PRK-DEIM}
satisfies
$$\norm{A(Nh) - Y_N}_F \leq c_{1,2} (\delta + c_{3,4} \cdot C_{m,r} C_{n,r} \cdot \varepsilon_r + h^q),$$
for every $0 < h < h_0$ and every $N \in \mathbb N$ such that $Nh \leq T$. 
The constants $c_{1,2}$, $c_{3,4}$ do not depend on $h$ or $\varepsilon_r$. 
\end{theorem}

\begin{proof}
The proof of \cite[Theorem 6]{kieri2019projection} holds verbatim, using Lemma~ \ref{lemma: projected stages error} for controlling the local stage error and $\varepsilon = C_{m,r} C_{n,r} \cdot \varepsilon_r$.
\end{proof}

The order $q$ established by Theorem~\ref{thm: convergence of PRK-DEIM} can be significantly smaller than the full order $p$ of the original RK method.
As discussed in~\cite{kieri2019projection}, the following PRK methods retain the full order $q = p = s$:
\[
\begin{aligned}
\text{projected Euler$=$PRK1:} &\quad b_1 = 1,\\[2pt]
\text{PRK2:} &\quad a_{21} = 1,\quad b_1 = b_2 = \tfrac{1}{2},\\[2pt]
\text{PRK3:} &\quad a_{21} = \tfrac{1}{3},\; a_{31} = 0,\; a_{32} = \tfrac{2}{3},\;
               b_1 = \tfrac{1}{4},\; b_2 = 0,\; b_3 = \tfrac{3}{4}.
\end{aligned}
\]


\subsection{Numerical experiment for a non-linear Schrödinger equation} \label{sec: schroedinger example}

All numerical experiments presented in the rest of this section have been carried out in Python on a MacBook Pro with an Apple M1 processor and 16GB RAM. More details, including the Python code with instructions on how to reproduce the experiments, are available on GitHub at \url{https://github.com/BenjaminCarrel/DLRA-DEIM}.

We consider the non-linear Schrödinger equation also considered in \cite{kieri2019projection} and derived from~\cite{trombettoni2001discrete}.
The solution $A: [0,T] \rightarrow \mathbb C^{n \times n}$ evolves according to the dynamics
\begin{equation} \label{eq: schrodinger}
\mathrm{i} \dt{A}(t) = - \frac{1}{2} (B A(t) + A(t) B) - \alpha A(t) * \overline{A(t)} * A(t),
\end{equation}
with $B= \mathrm{tridiag}(1, 0, 1)$. The initial value is given by the rank-$2$ matrix
\begin{equation*}
A_{j,k}(0) = \exp \left(- \frac{(j- \mu_1)^2}{\sigma^2} - \frac{(k-\nu_1)^2}{\sigma^2} \right) + \exp \left( - \frac{(j-\mu_2)^2}{\sigma^2} - \frac{(k - \nu_2)^2}{\sigma^2} \right),
\end{equation*}
with scaled parameters $\sigma = 0.1n, \mu_1 = 0.6n, \mu_2 = 0.5n, \nu_1 = 0.5n,$ and $\nu_2 = 0.4n$.
As in \cite{kieri2019projection}, we first propagate the solution from $t=0$ to $t=0.01$ to ensure full rank for the initial value. The final time is $T=1$ and we set $\alpha = 0.1$, $n=1024$.

Because~\eqref{eq: schrodinger}  is complex-valued, we need to apply complex versions of the DEIM procedures. While QDEIM and Osinsky’s method are originally given for complex matrices, ARP can be extended in a straightforward way. The bounds stated in Table~\ref{tab:algo_comparison} and the error analysis of PRK-DEIM hold for the complex extensions of these methods as well.  In the experiments, we used ARP since it gave the smallest error with Osinsky’s method a close second at essentially the same computational time.

\begin{figure}[!ht]
\centering
\includegraphics[width=0.85\textwidth]{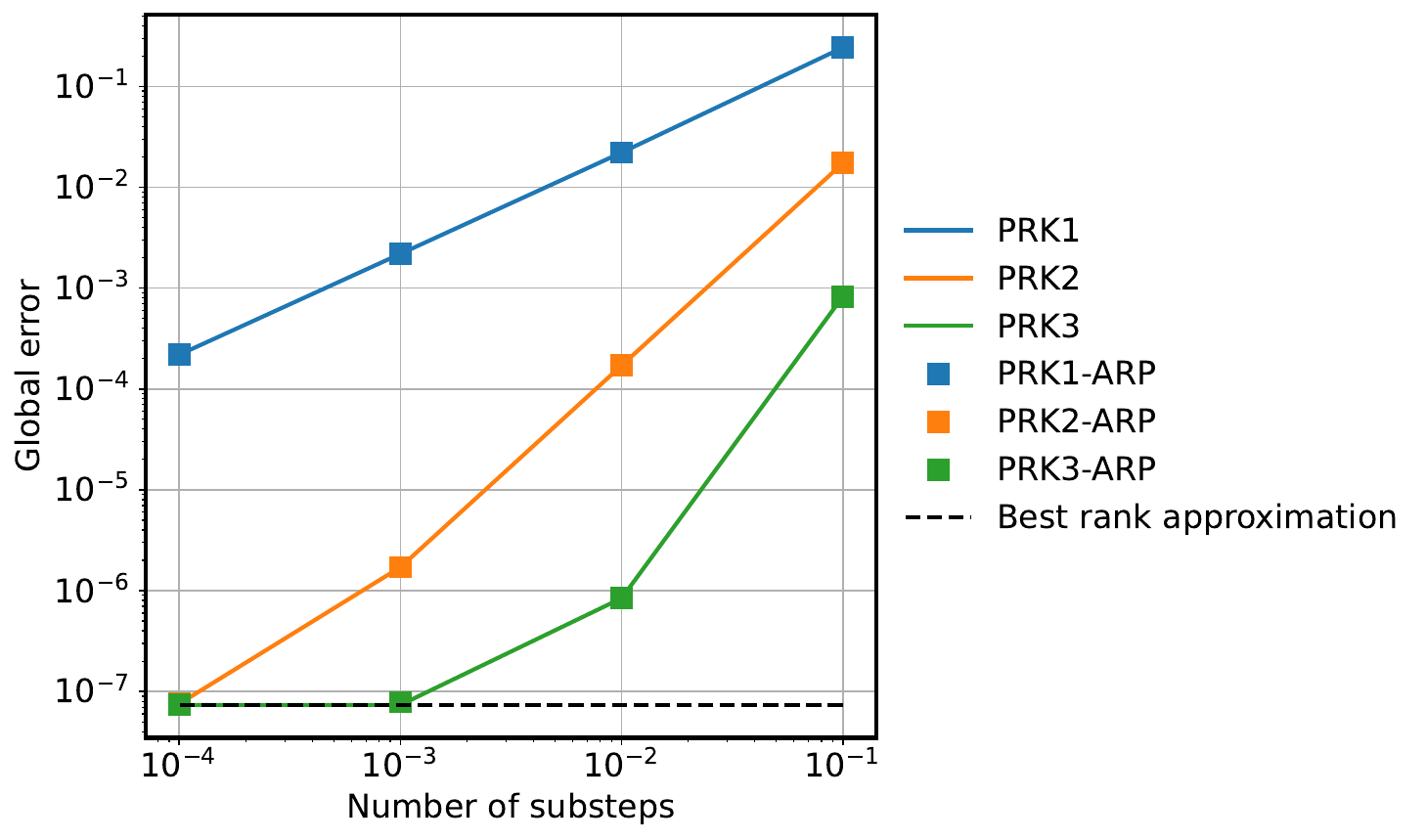}
\caption{Relative error at $T = 1$ of various PRK and PRK-DEIM methods with ARP applied to the Schrödinger equation \eqref{eq: schrodinger} with approximation rank $r=8$. }
\label{fig: PRK-DEIM global error}
\end{figure}

Figure \ref{fig: PRK-DEIM global error} confirms our theory; the new PRK-DEIM methods exhibit the same order as the PRK methods, up to the modeling error.
In fact,
the errors of PRK-DEIM visibly match the errors of PRK. While the 
theory does not exclude a moderate increase of the (modeling) error due to oblique DEIM projections, such an increase is not observed in the figure. Here and in the following, the relative error is the Frobenius norm of the difference between the reference solution and the approximated solution, normalized by the Frobenius norm of the reference solution.

Table~\ref{table: performance PRK} provides a more detailed picture, comparing  accuracy and efficiency. When executing the original PRK methods, we evaluated the non-linearity $A * \overline{A} * A$ for a rank-$r$ matrix $A$ efficiently by representing this rank-$r^3$ matrix in factorized form, as described in \cite{kressner2017recompression}. However, already for relatively small ranks $r$ the increase to $r^3$ is costly. On the other hand, evaluating $A * \overline{A} * A$ element-wise at the indices selected by the DEIM procedure is much cheaper. This is clearly reflected in the timings of Table~\ref{table: performance PRK}; with a more pronounced effect as $r$ gets larger.

\begin{table}[h]
\centering
\begin{tabular}{l rrr} 
\toprule
Method & $r=3$ & $r=6$ & $r=9$ \\
\midrule
\multicolumn{4}{l}{\textit{Relative errors at final time}} \\ \midrule
PRK1         & 7.8666e-03 & 2.1883e-03 & 2.1882e-03 \\ 
PRK2         & 7.5486e-03 & 2.6146e-05 & 1.7120e-06 \\ 
PRK3         & 7.5486e-03 & 2.6090e-05 & 7.3686e-08 \\
PRK1-ARP     & 7.9453e-03 & 2.1880e-03 & 2.1882e-03 \\
PRK2-ARP     & 7.5657e-03 & 2.6554e-05 & 1.7110e-06 \\
PRK3-ARP     & 7.5700e-03 & 2.6720e-05 & 7.6915e-08 \\
\midrule 
\multicolumn{4}{l}{\textit{Execution time (in seconds)}} \\ \midrule 
PRK1         & 12.12 & 65.92 & 126.23 \\ 
PRK2         & 24.11 & 133.50 & 257.43 \\ 
PRK3         & 38.36 &  205.52 &  390.94 \\
PRK1-ARP     & 2.62  &   5.54 &   9.86 \\
PRK2-ARP     & 5.88  &   14.62 &  24.12 \\
PRK3-ARP     & 11.00 &  25.19 &   44.53 \\
\bottomrule
\end{tabular}
\caption{Performance comparison of various PRK and PRK-DEIM methods with ARP, applied to the Schrödinger equation \eqref{eq: schrodinger} with step size $h=10^{-3}$ and final time $T = 1$.}
\label{table: performance PRK}
\end{table}

\subsection{Exponential integrators}

PRK and PRK-DEIM are both explicit time stepping methods that will suffer from step size restrictions when applied to stiff problems, as they arise from, e.g., the discretization of time-dependent partial differential equations. One way to address this issue is to use implicit DLRA schemes~\cite{appelo2024robust,rodgers2023implicit,nakao2023reduced}.
When the vector field $F$ admits a splitting of the form $F=L+G$, with a linear part $L$ that accounts for most of the stiffness and a non-linear part $G$, 
there are appealing alternatives to implicit methods. In particular, when $L$ is a so-called Sylvester operator (e.g., a discretized 2D-Laplace operator) then 
explicit exponential Runge--Kutta methods combine nicely with low-rank structure~\cite{ostermann2019convergence,carrel2023low}. Because they are close in spirit to PRK, we consider the so-called projected exponential Runge--Kutta (PERK) methods
from~\cite{carrel2023projected}, which combine the explicit exponential Runge--Kutta methods from~\cite{Hochbruck2005} with rank-$r$ truncation and orthogonal tangent space projections. The convergence of PERK is provably not affected by stiffness of $L$.

PERK-DEIM methods are obtained from PERK methods by replacing orthogonal with DEIM-based oblique tangent space projections, just as PRK-DEIM methods were obtained from PRK methods in~\eqref{eq: PRK-DEIM}. This construction is straightforward and, therefore, we do not provide details. In complete analogy to the analysis in Section~\ref{sec:error}, one can show that the convergence results from~\cite{carrel2023projected} hold with minor modifications. In particular, 
the order of convergence is not affected when using DEIM.

\subsection{Numerical experiment for Allen--Cahn equation}
\label{sec: Allen-cahn_matrix}
To numerically validate PERK-DEIM, we now consider the Allen--Cahn equation from \cite{rodgers2023implicit}, simulating phase separation in multi-component alloy systems \cite{allen1972ground, allen1973correction}. 
In its simplest form, this equation is given by
\begin{equation}
\label{eq ac-eq}
\frac{\partial f}{\partial t} = \kappa \Delta f + f - f^3,
\end{equation}
where $\Delta f$ is the diffusion term (with the Laplace operator $\Delta$), $f - f^3$ is the reaction term, and $\kappa \in \R$ us a parameter. This (discretized) diffusion term makes this a stiff problem.

We consider the spatial domain $\Omega = [0, 2 \pi]^2$ with periodic boundary conditions, and the time interval $[0, T] = [0, 10]$.
The spatial discretization is done via central finite differences, leading to a matrix differential equation of the form
\begin{equation} \label{eq: Allen-Cahn}
\begin{aligned}
&\dt{A}(t) = D A(t) + A(t) D + A(t) - A(t)*A(t)*A(t), \\
&A(0)_{ij} = \frac{\left[ e^{-\tan^2(x_i)} + e^{-\tan^2(y_j)} \right] \sin(x_i) \sin(y_j)}{1 + e^{|\csc(-x_i/2)|} + e^{|\csc(-y_j/2)|}}
\end{aligned}
\end{equation}
where the tridiagonal matrix $D \in \Rnn$ is the discretized one-dimensional Laplace operator with periodic boundary conditions, scaled by $\kappa$. As explained in Section~\ref{sec: schroedinger example}, the evaluation of the cubic nonlinearity 
$A*A*A$ for a rank-$r$ matrix $A$ greatly benefits from the use of DEIM.

To assess performance, we have considered $n=256$ and $n=1024$ with parameter $\kappa = 0.01$. These choices result in stiff problems that lead to unreasonable step size restriction for standard explicit methods. In contrast, we expect that PERK and PERK-DEIM are not subject to these restrictions. Indeed, the results 
of Table~\ref{table: performance PERK} show that various instances of these methods result in good accuracy for the relatively large step size $h = 10^{-2}$.
Compared to the methods from~\cite{carrel2023projected} that utilize orthogonal projections (PERK1 and PERK2), our new methods with oblique projections are significantly faster, with the difference becoming more dramatic as $r$ increases. 
We have included results for two DEIM procedures, SRRQR and ARP; the resulting cost and accuracy are comparable. 
Note that the time integration of the linear part in 
PERK and PERK-DEIM requires one to apply the matrix exponential of $h D$ to the low-rank factors of the current iterate. 
We have used an extended Krylov subspace method of dimension two for this purpose, ensuring that the error due to the Krylov subspace approximation is negligible compared to the modeling error.

\begin{table}[h]
\centering
\begin{tabular}{l rrrr} 
\toprule
Method & \multicolumn{2}{c}{$n=256$} & \multicolumn{2}{c}{$n=1024$} \\
\cmidrule(lr){2-3} \cmidrule(lr){4-5} 
       & $r=3$ & $r=6$ & $r=6$ & $r=9$ \\
\midrule
\multicolumn{5}{l}{\textit{Relative errors at final time}} \\ \midrule
PERK1         & 2.3228e-03 & 5.4594e-04 & 5.4580e-04 & 4.5209e-04 \\ 
PERK2         & 2.2784e-03 & 3.0579e-04 & 3.0554e-04 & 5.5290e-06 \\ 
PERK1-SRRQR  & 2.3048e-03 & 6.2657e-04 & 6.2730e-04 & 4.5309e-04 \\ 
PERK2-SRRQR  & 2.2754e-03 & 3.7759e-04 & 3.7788e-04 & 1.7211e-05 \\ 
PERK1-ARP     & 2.3201e-03 & 4.8066e-04 & 4.6852e-04 & 4.5211e-04 \\
PERK2-ARP     & 2.2786e-03 & 1.2126e-04 & 1.5431e-04 & 5.8688e-06 \\ 
\midrule 
\multicolumn{5}{l}{\textit{Execution time (in seconds)}} \\ \midrule 
PERK1         &  1.94 &   8.31 &  32.28 &  110.75 \\ 
PERK2         & 4.87 &  20.31 &  75.67 &  225.09 \\ 
PERK1-SRRQR  &  1.79 &   3.01 &  14.60 &   17.37 \\
PERK2-SRRQR  & 4.63 &   9.13 &  39.03 &   52.01 \\ 
PERK1-ARP     & 1.74 &   2.85 &  10.51 &   17.76 \\ 
PERK2-ARP     & 4.45 &   9.18 &  38.84 &   53.60 \\
\bottomrule
\end{tabular}
\caption{Performance comparison of various PERK and PERK-DEIM methods applied to the Allen--Cahn equation \eqref{eq: Allen-Cahn} with step size $h=10^{-2}$.}
\label{table: performance PERK}
\end{table}

\section{Extension to tensor differential equations} \label{sec:tucker}

The purpose of this section is to show how PRK-DEIM methods extend from matrix to low-order tensor differential equations by adapting the DEIM-based oblique tangent space projections.

To simplify the description, we consider third-order tensors with identical mode sizes: $A(t)\in\R^{n \times n \times n}$. The extension to arbitrarily-sized tensors is straightforward.
Consider the differential equation \begin{equation}
\label{eq: tensor ode 1}
    \dot{A}(t)=F(A(t)),\quad A(0)=A_0\in \R^{n \times n \times n}.
\end{equation} We aim at approximating its solution $A(t)$ by a tensor $Y(t)\in\R^{n \times n \times n}$ of low multilinear rank. Again for simplicity, we assume that the desired multilinear rank is identical to the same value $r\ll n$ for each mode. In \cite{Koch2010dynamical}, DLRA was extended to tensors by considering the manifold 
$$\M_r=\{Y\in \R^{n\times n\times n}\,|\, Y\text{ has multilinear rank $r$}\},$$ to  obtain $Y(t)$ from solving the differential equation
\begin{equation}\label{eq: tensor ode}\dot{Y}(t)=\mathcal{P}_{Y(t)}[F(Y(t))],\quad Y(0)=Y_0,\end{equation} where $Y_0$ is a multilinear rank-$r$ approximation to $A_0$ and $\mathcal{P}_{Y(t)}[\cdot]$ is the orthogonal projection onto the tangent space $T_{Y(t)}\mathcal{M}_r$.

To turn~\eqref{eq: tensor ode} into an efficient method, one leverages the fact that a tensor $Y$ of low multilinear rank can be compactly represented in the so-called Tucker format. In the following, we briefly recall the necessary ingredients to define this format and refer to \cite[Part II]{Ballard2025tensor} for details.
The mode-$i$ unfolding of a general tensor $B\in \R^{n_1\times n_2\times \cdots \times n_d}$ is the matrix $B_{(i)}\in \R^{n_i\times (n_1\cdots n_{i-1} n_{i+1}\cdots n_d)}$ obtained from stacking the mode-$n_i$ fibers of $B$. Given a matrix 
$V\in \R^{m\times n_i}$, this allows us to express 
the mode-$i$ product $B\times_i V \in \R^{n_1\times \cdots\times  n_{i-1}\times m\times n_{i+1}\times \cdots \times n_d}$ in terms of a matrix-matrix product: $$(B\times_i V)_{(i)}=V B_{(i)}.$$
A tensor $Y \in \R^{n\times n\times n}$ has multilinear rank $r$ if each of its unfoldings $Y_{(1)}$, $Y_{(2)}$, $Y_{(3)}$ has rank $r$. Moreover, such a tensor
admits the following representation in the Tucker format: 
\begin{equation}
\label{eq: tucker decomposition}
Y =C \times_1U_1 \times_2U_2 \times_3U_3 =: C \times^3_{i=1}U_i,
\end{equation} with orthonormal $U_i(t)\in \R^{n\times r}$ and the so called core tensor $C \in\R^{r\times r\times r}$. This format can be turned into
rank-$r$ factorizations for the unfoldings as follows:
\begin{equation} \label{eq:factucker}
    Y_{(i)} = U_i S_i V_i^\top,
\end{equation}
with the right factors defined by the orthonormal $n^2\times r$ matrices
\begin{equation} \label{eq:Vi}
 V_1 = (U_3\otimes U_2)\, Q_1, \quad V_2 = (U_3\otimes U_1)\, Q_2, \quad V_3 = (U_2\otimes U_1)\, Q_3,
\end{equation}
where ``$\otimes$'' denotes the Kronecker product.
The matrices $Q_i \in \mathbb{R}^{r^2 \times r}$
and $S_i \in \mathbb{R}^{r \times r}$ are obtained from an economy-sized QR factorization of the mode-$i$ unfolding of the core tensor: $C_{(i)}^\top = Q_i S_i^\top$. 

Following~\cite[Section 6]{lubich2015time}, we use~\eqref{eq:factucker} to express the orthogonal tangent space projection $\mathcal{P}_{Y}$ as a sum of partial projections. For this purpose, we let $\text{ten}_{i}$ denote the inverse of mode-$i$ unfolding, which is uniquely defined by the relation $B=\text{ten}_i(B_{(i)})$ for an arbitrary $n\times n\times n$ tensor $B$. Defining
\begin{equation}\label{eq: subprojection 1}
    P_{i}(Y) [Z] = \text{ten}_i \big( (I - U_i U_i^\top) Z_{(i)} {V}_i V_i^\top\big), \quad
    \bar P(Y) [Z] = Z \times^3_{i=1} U_i U_i^\top,
\end{equation}
we have that
$\mathcal{P}_{Y}[Z]=P_{1}(Y)[Z] + \cdots + P_{3}(Y)[Z] + \bar P(Y) [Z].$

\subsection{Oblique tangent space projections for Tucker tensor}

The computation of the partial projections~\eqref{eq: subprojection 1} requires full knowledge of the tensor $Z$. Like as in the matrix case, this may incur a significant cost when the differential equation~\eqref{eq: tensor ode 1} has nonlinearities. This can be avoided by replacing the orthogonal projections $U_i U_i^\top$,
$V_i V_i^\top$ appearing in~\eqref{eq: subprojection 1} by the 
DEIM-based oblique projections $P_{U_i}^{\angle}$, $P_{V_i}^{\reflectangle}$ defined in \eqref{eq: oblique projections range and corange}:
\begin{equation}
{P}^{\angle}_{i}(Y) [Z] := \text{ten}_i \big( (I-{P}^{\angle}_{U_i}) Z_{(i)} P_{V_i}^{\reflectangle}\big),\quad \bar{P}^{\angle}(Y) [Z] := Z \times^3_{i=1} {P}^{\angle}_{U_i}.
\end{equation}
In turn, the orthogonal projection $\mathcal{P}_{Y}[Z]$ is replaced by 
\begin{equation}
\label{eq: tensor oblique projection}
    \mathcal{P}^{\angle}_Y[Z]=
    P^{\angle}_{1}(Y)[Z] + \cdots + P^{\angle}_{3}(Y)[Z] + \bar{P}^{\angle}(Y)[Z].
\end{equation}
It is not difficult to see that ${O}(rn+r^3)$ entries of $Z$ need to be evaluated to define $\mathcal{P}^{\angle}_Y[Z]$: $r$ mode-$i$ fibers for each ${P}^{\angle}_{i}(Y) [Z]$ and $r^3$ entries for $\bar{P}^{\angle}(Y)$.

The following lemma lists basic properties of $\mathcal{P}^\angle_{Y}[Z]$. Since $\mathcal{P}^\angle_{Y}[Z]$ is linear in $Z$, they imply that $\mathcal{P}^\angle_{Y}$ is an oblique projector onto $T_{Y}\mathcal{M}_r$.
\begin{lemma}\label{lemma:tensorproj}Assume that $P_{U_i}^{\angle}$ and $P_{V_i}^{\reflectangle}$ are well defined, that is, $S^\top_{U_i}U_i$ and $S^\top_{V_i}V_i$ are invertible for $i=1,2,3$. Then
   \begin{enumerate}
       \item  $\mathcal{P}^{\angle}_Y[\mathcal{P}^{\angle}_Y[Z]]=\mathcal{P}^{\angle}_Y[Z]$,
       \item $\mathcal{P}_Y[\mathcal{P}^{\angle}_Y[Z]]=\mathcal{P}^{\angle}_Y[Z]$,
       \item $\mathcal{P}^{\angle}_Y[\mathcal{P}_Y[Z]]=\mathcal{P}_Y[Z]$.
   \end{enumerate}
\end{lemma}
\begin{proof} We only prove point 3.; the proofs of 1. and 2. follow in an analogous fashion. We first note that the definition~\eqref{eq:Vi} of $V_j$ implies that
the mode-$i$ fibers of ${P}_{j}(Y) [Z]$ are contained in the range of $U_i$ for $j \not= i$. Because $P_{U_i}^{\angle}$ is a projector onto that range, it follows that
\[
 P^{\angle}_{i}(Y)[ P_{j}(Y)[Z] ] = 0, \quad \text{for $j\not= i$}, \quad 
 P^{\angle}_{i}(Y)[ \bar P(Y)[Z] ] = 0.
\]
This implies
\begin{align*}
 P^{\angle}_{i}(Y)[ \mathcal{P}_Y[Z] ] &=  P^{\angle}_{i}(Y)[ P_{i}(Y)[Z] ] 
 =
 \text{ten}_i \big( (I-{P}^{\angle}_{U_i}) (I-U_i U_i^\top) Z_{(i)} V_i V_i^\top P_{V_i}^{\reflectangle} \big) \\
 &= \text{ten}_i \big( (I-{P}^{\angle}_{U_i}) Z_{(i)} V_i V_i^\top  \big),
\end{align*}
where we used~\eqref{eq:projectors} for the last equality. Similarly, we obtain that
\[
 \bar{P}^{\angle}(Y) [ \bar{P}(Y)[Z] ] = \bar{P}(Y)[Z], \quad 
 \bar{P}^{\angle}(Y) [ P_{i}(Y)[Z] ] = \text{ten}_i \big( ( {P}^{\angle}_{U_i} - U_i U_i^\top) Z_{(i)} V_i V_i^\top \big).
\]
Taking all relations into account, this gives 
\[
 \mathcal{P}^{\angle}_Y[\mathcal{P}_Y[Z]] = 
 \text{ten}_i \Big(
 \sum_{i = 1}^3
 (I-{P}^{\angle}_{U_i}) Z_{(i)} V_i V_i^\top  + 
 ( {P}^{\angle}_{U_i} - U_i U_i^\top) Z_{(i)} V_i V_i^\top
 \Big) + \bar{P}(Y)[Z] = \mathcal{P}_Y[Z],
\]
as desired.
\end{proof}

\subsubsection{Index selection for $V_i$}

Our construction~\eqref{eq: tensor oblique projection} of $\mathcal{P}^{\angle}_Y$ involves three oblique projectors of the form $P_{V}^{\reflectangle} = S_{V} (V^\top S_{V})^{-1} V^\top$, with an $n^2\times r^2$ matrix
$V=(\bar{U}\otimes \tilde{U})Q$ for orthonormal matrices $\bar{U}, \tilde{U} \in \mathbb{R}^{n \times r}$ and $Q \in \mathbb{R}^{r^2 \times r}$.
Applying a DEIM procedure to $V$ directly is expensive; for example, QDEIM requires $O(n^2r^4)$ operations. Instead, we suggest to apply
DEIM in two stages: First to $\bar{U}$ and $\tilde{U}$ to obtain selection matrices $S_1$ and $S_2$, respectively, and then to an orthonormal basis of $[(S_1^\top \bar{U})\otimes (S_2^\top \tilde{U})]Q\in \mathbb{R}^{r^2\times r}$, yielding a selection matrix $S_{1,2}\in \mathbb{R}^{r^2\times r}$. This gives the final selection matrix $S_V=(S_1\otimes S_2)S_{1,2}$. 
\begin{algorithm}[H]
\caption{Construction of sampling matrix \(S_V\) for $V=(\bar{U}\otimes \tilde{U})Q$}
\label{alg: sampling v}
\begin{algorithmic}[1]
\Require Orthonormal matrices \(\bar{U} \in \mathbb{R}^{n \times r}\), \(\tilde{U}\in \mathbb{R}^{n \times r}\), and \(Q \in \mathbb{R}^{r^2 \times r}\)
\State \(S_1 \gets \text{DEIM}(\bar{U}),\quad S_2 \gets \text{DEIM}(\tilde{U})\) 
\State Compute orthonormal basis \(\hat{Q} \gets \mathsf{orth}([(S_1^\top \bar{U}) \otimes (S_2^\top \tilde{U})]Q)\) by QR decomposition
\State \(S_{1,2} \gets \text{DEIM}(\hat{Q})\)
\State Return \(S_V \gets (S_1 \otimes S_2)\,S_{1,2}\) \Comment{Stored implicitly via row indices of $S_1, S_2, S_{1,2}$.}
\end{algorithmic}
\end{algorithm}%
Algorithm~\ref{alg: sampling v} summarizes the described procedure. When using, e.g.,
QDEIM, Osinsky's method, or ARP for DEIM, this algorithm requires $O(nr^2 + r^4)$ operations.
Lemma~\eqref{lemma: algo2} below establishes a bound on $\|P_{V}^{\reflectangle}\|_2$, measuring the quality of the oblique projection. 
\begin{lemma}
\label{lemma: algo2}
    The sampling matrix $S_V$ returned by Algorithm \ref{alg: sampling v} yields an oblique projector $P_{V}^{\reflectangle} = S_{V} (V^\top S_{V})^{-1} V^\top$
    such that 
    $\|P_{V}^{\reflectangle}\|_2 \leq C^2_{n,r}C_{r^2,r}$ with $C_{m,r}$ from~\eqref{eq:def Cmr}.
\end{lemma}
\begin{proof}
    Considering the QR decomposition $\hat{Q}\hat{R}=[(S_1^\top \bar{U})\otimes (S_2^\top \tilde{U})]Q$ and using basic properties of the Kronecker product, we obtain that
    \begin{align*}
     \|P_{V}^{\reflectangle}\|_2 &= \|(S^\top V)^{-1}\|_2 = \|\big(S_3^\top[(S_1^\top \bar{U})\otimes (S_2^\top \tilde{U})]Q\big)^{-1}\|_2 \\
    &= \|\hat{R}^{-1}(S_3^\top \hat{Q})^{-1}\|_2
 \leq \|[(S_2^\top \bar{U})\otimes (S_2^\top \tilde{U})]^{-1}\|_2\|(S_3^\top \hat{Q})^{-1}\|_2\\&= \|(S_1^\top \bar{U})^{-1}\|_2\|(S_2^\top \tilde{U})^{-1}\|_2\|(S_3^\top \hat{Q})^{-1}\|_2
 \leq C^2_{n,r}C_{r^2,r}.
 \end{align*}
\end{proof}
For example, for Osinsky's method or ARP, we obtain from Table~\ref{tab:algo_comparison} that (in expectation) 
\[
 \|P_{V}^{\reflectangle}\|_2 \leq  (1 + r (n-r)) \sqrt{ 1 + r^2 (r^2-r) } \sim n r^3.
\]

\subsection{Error analysis of the oblique tangent space projection}

The following lemma provides a quasi-optimality result for $P^{\angle}_Y[Z]$.
Note that $\|Z\|_F$ denotes the usual Frobenius for a tensor
$Z\in \R^{n\times n \times n}$: $\|Z\|_F^2 = \sum_{ijk} Z_{ijk}^2$.
\begin{lemma}
\label{lemma: Tucker projection error}
Suppose that $S^\top_{U_i}U_i$ and $S^\top_{V_i}V_i$ are invertible for all $i=1,2,3$. Then the oblique projector $\mathcal P^{\angle}_Y[Z]$ defined in~\eqref{eq: tensor oblique projection} satisfies
    \begin{align}\|\mathcal P^{\angle}_Y[Z]-Z\|_F\leq\Big(\sum_{i = 1}^3 \| P^{\angle}_{U_i}\|_2 \| P^{\reflectangle}_{V_i}\|_2 + \prod_{i = 1}^3 \| P^{\angle}_{U_i}\|_2\Big)\|\mathcal  P_Y[Z]-Z\|_F 
\end{align}
In particular, when using Algorithm \ref{alg: sampling v} for selecting the rows of $V_i$, it holds that
\begin{equation}
\label{eq: constant by algorithm 2}
    \|\mathcal P^{\angle}_Y[Z]-Z\|_F\leq
     ( 3 C_{r^2,r} + 1) C^3_{n,r} \|\mathcal P_Y[Z]-Z\|_F.
\end{equation}
\end{lemma}
\begin{proof}
By Lemma~\ref{lemma:tensorproj}, $\mathcal P^{\angle}_Y[\mathcal  P_Y[Z]]=\mathcal P_Y[Z]$ and thus
\[
 \|\mathcal P^{\angle}_Y[Z]-Z\|_F = \|(\mathcal{I}-\mathcal P^{\angle}_Y)[\mathcal P_{Y}[Z]-Z]\|_F \le \|\mathcal{I}-\mathcal P^{\angle}_Y\|_2\|\mathcal{P}_Y[Z]-Z\|_F = 
 \|\mathcal P^{\angle}_Y\|_2\|\mathcal{P}_Y[Z]-Z\|_F,
\]
where $\|\cdot\|_2$ denotes the operator norm induced by $\|\cdot\|_F$. Using the definition~\eqref{eq: tensor oblique projection},
we have
\[
 \|\mathcal P^{\angle}_Y\|_2 \le \sum_{i = 1}^3 \| P^{\angle}_{i}(Y) \|_2 + \|
 \bar P^{\angle}(Y) \|_2 \le \sum_{i = 1}^3 \| P^{\angle}_{U_i}\|_2 \| P^{\reflectangle}_{V_i}\|_2 + \prod_{i = 1}^3 \| P^{\angle}_{U_i}\|_2, 
\]
which shows the first part of the lemma. The inequality~\eqref{eq: constant by algorithm 2} is obtained by applying Lemma~\ref{lemma: algo2}.
\end{proof}

\subsection{Tucker PRK-DEIM}

Tucker PRK-DEIM is a direct generalization of PRK-DEIM~\eqref{eq: PRK-DEIM}
to the 
tensor differential equation~\eqref{eq: tensor ode 1}. The only notable differences are that the truncation to the manifold $\mathcal M_r$ of third-order tensors having multilinear rank $r$ is performed by applying the 
quasi-optimal higher-order SVD (HOSVD) from~\cite{De2000multilinear}, and the oblique tangent space projections are carried out by applying~\eqref{eq: tensor oblique projection}, using a DEIM procedure combined with Algorithm~\ref{alg: sampling v} for index selection. The error analysis of Tucker PRK-DEIM follows from a straightforward extension of Theorem~\ref{thm: convergence of PRK-DEIM}, using the constants from Lemma~\ref{lemma: Tucker projection error}.

 \subsection{Numerical experiments}
 
  In this section, we report two numerical experiments for various instances of Tucker PRK-DEIM. We computed reference solutions by using the classical fourth-order Runge–Kutta method with a small time step, $h = 10^{-3}$. All experiments in this section have been performed in Matlab (version 2023a) on a Macbook Pro with an
Apple M1 Pro processor. The code used to perform the experiments and reproduce the figures in this section can be found at \url{https://github.com/hysanlam/DLRA_DEIM_tensor}.
 \subsubsection{A discrete nonlinear Schrödinger equation}

\begin{figure}[H]
     \centering
     \includegraphics[scale=0.65]{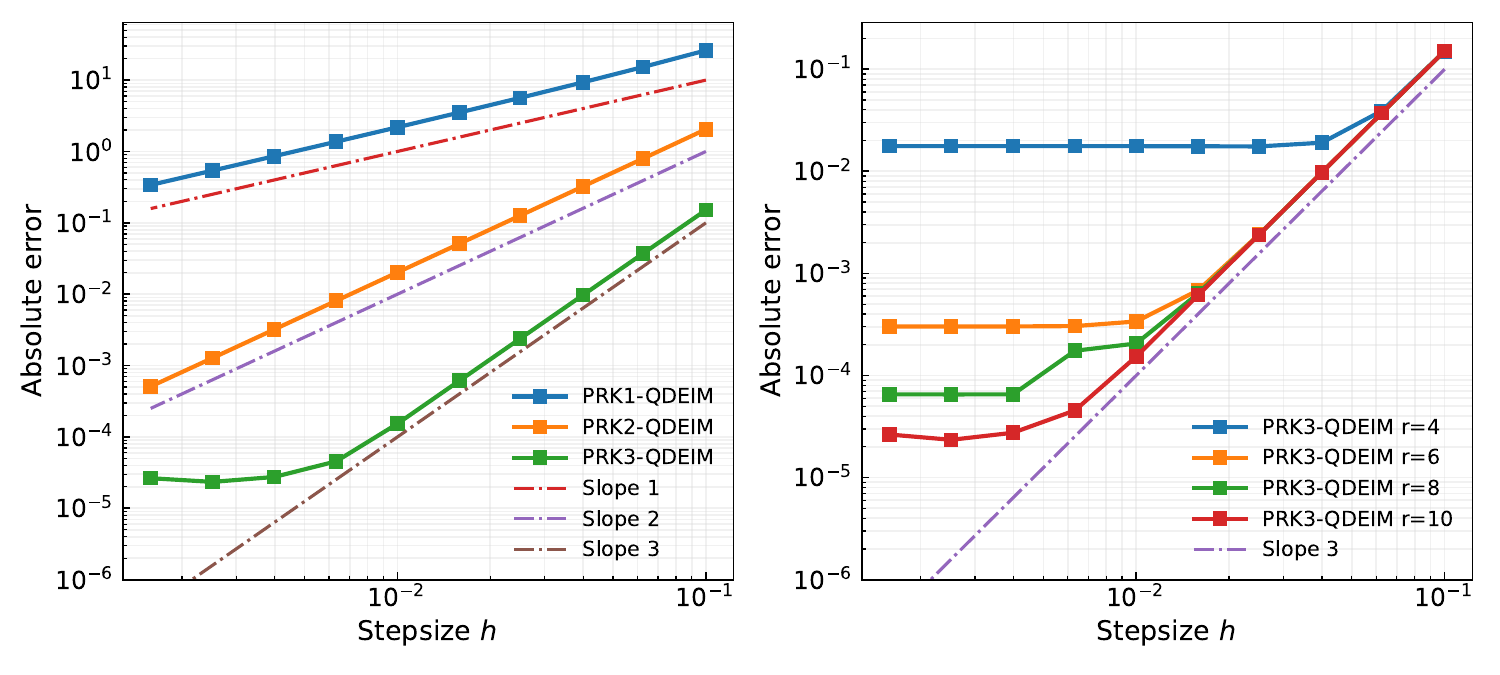}
     \caption{Discrete nonlinear Schrödinger equation~\eqref{eq: tensor schrodinger}: Frobenius norm error of Tucker PRK-QDEIM  methods vs time-step size at $T=1.01$ for fixed rank (left panel) and varying rank (right panel). }
     \label{fig:tensor schrodinger}
 \end{figure}

We use an example from \cite{Lubich2018time}, which is a three-dimensional variation of the example considered in Section~\ref{sec: schroedinger example}, and consider the discrete nonlinear Schrödinger equation
\begin{equation}
\label{eq: tensor schrodinger}
    \mathrm{i} \dot{A}(t) = -\frac{1}{2}L[A(t)] + \alpha A(t) * \overline{A(t)} * A(t),
\end{equation}
with the initial value
\begin{align*}
    A_{jkl}(0) &= \exp \left(-1/\gamma^2 \left( (j - 75)^2 - (k - 25)^2 - (l - 1)^2 \right) \right) 
    \\&+ \exp \left(-1/\gamma^2 \left( (j - 25)^2 - (k - 75)^2 - (l - 100)^2 \right) \right).
\end{align*}
We choose a lattice size $n = 100$ for each coordinate, such that
$A(t) \in \mathbb{R}^{n \times n \times n}$, and $\alpha=10^{-1}$. The linear operator 
$L: \mathbb{R}^{n \times n \times n} \to \mathbb{R}^{n \times n \times n}$ describes the interaction between the grid points centered at $(j,k,l)$ for all $j,k,l = 1, \dots, 100$. It is defined component-wise as
\begin{equation}
    L[A](j,k,l) = A(j-1,k,l) + A(j+1,k,l) + A(j,k-1,l) + A(j,k+1,l)
    + A(j,k,l-1) + A(j,k,l+1).
\end{equation}
The nonlinearity $A(t) * \overline{A(t)} * A(t)$ in~\eqref{eq: tensor schrodinger} is evaluated entry-wise. Again, 
we first propagate the solution from $t=0$ to $t=0.01$ using the classical fourth-order Runge–Kutta method to ensure full multilinear rank $r$ for the (truncated) initial value. We have combined Tucker PRK-DEIM with the three Runge--Kutta methods PRK1, PRK2, PRK3 from Section~\ref{sec:error} and QDEIM. 
In the left panel of Figure \ref{fig:tensor schrodinger}, we provide the absolute error for fixed multilinear rank $r=10$ while varying the time-step size. One observes that all methods exhibit the order of convergence (up to modeling error) predicted by Theorem \ref{thm: convergence of PRK-DEIM}. The right panel shows the behavior of PRK3-QDEIM  for different multilinear ranks. Again, one observes convergence at the right rate until the low-rank approximation error dominates.
 
 \subsubsection{3D Allen--Cahn equation}
 
This experiment uses the Allen--Cahn equation~\eqref{eq ac-eq} for the three-dimensional domain $\Omega = [0,2\pi]^3$, $\kappa = 0.1$, and the initial condition 
\[
f_0(x_1, x_2, x_3) = g(x_1, x_2, x_3) - g(2x_1, x_2, x_3) + g(x_1, 2x_2, x_3) - g(x_1, x_2, 2x_3),
\]
where
\[
g(x_1, x_2, x_3) = \frac{\left(e^{-\tan(x_1)^2} + e^{-\tan(x_2)^2} + e^{-\tan(x_3)^2} \right) \sin(x_1 + x_2 + x_3)}
{1 + e^{|\csc(-x_1/2)|} + e^{|\csc(-x_2/2)|} + e^{|\csc(-x_3/2)|}}.
\]
We discretize $\Omega$ using $n = 150$ points in each coordinate and approximate the Laplace operator using second-order finite differences. To avoid instability, we again integrate the tensor differential equation with the classical fourth-order Runge--Kutta method until $t=1$ and apply our method from there.

\begin{figure}[H]
     \centering
     \includegraphics[scale=0.65]{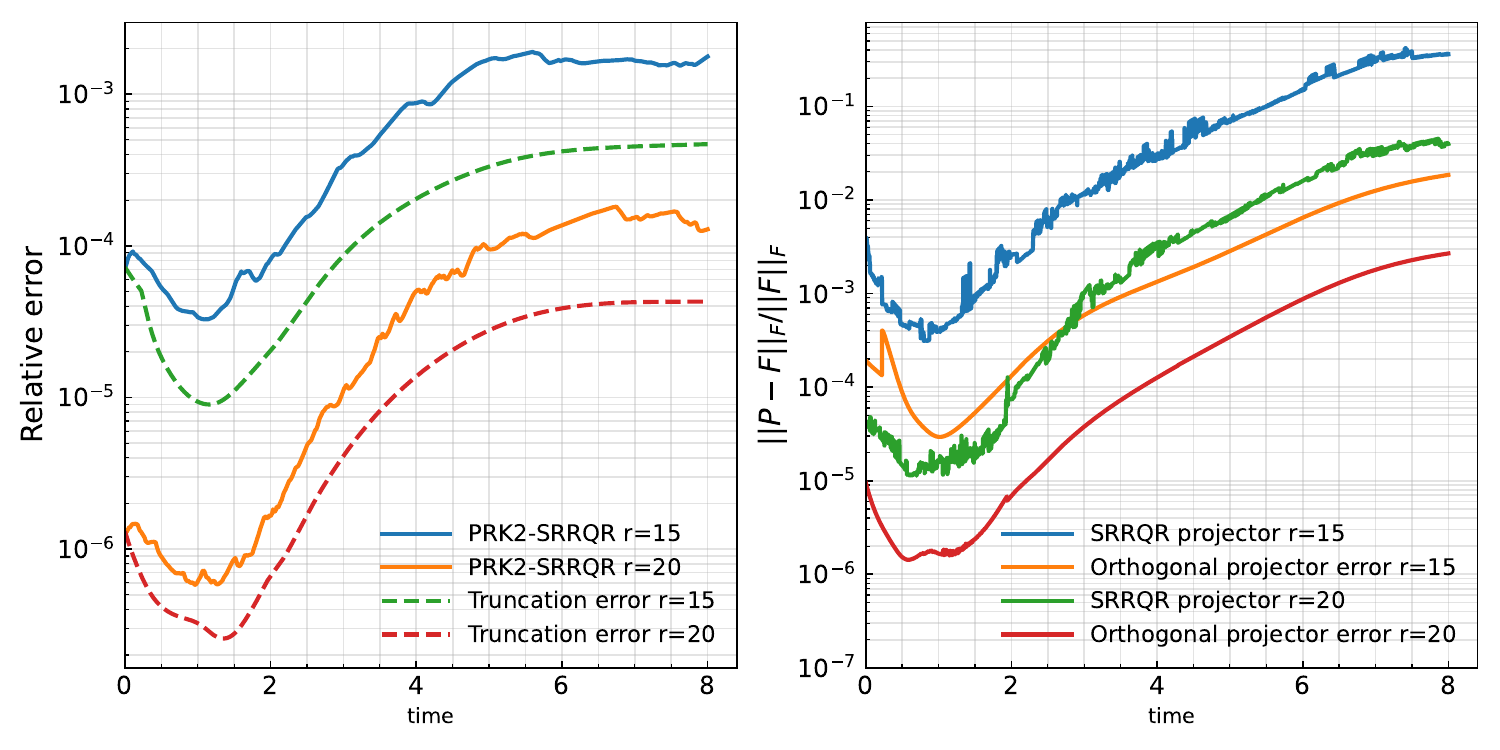}
     \caption{Results for 3D Allen--Cahn equation. Left: Relative error for Tucker PRK2-SRRQR $(\eta=2, h=10^{-3})$ error and best low rank approximation of reference solution. Right: Relative error of oblique and orthogonal projections of the vector field at the HOSVD truncation of the reference solution.}
     \label{fig:ac-3d}
 \end{figure}

The left panel of Figure~\ref{fig:ac-3d} shows the relative error of PRK2-SRRQR with multilinear rank $r = 15$ and $r = 20$, with respect to the shifted time $t-1$, compared to the (full rank) reference solution.
In the same plot, we also see the quasi-best truncation relative error, that is, the relative error obtained when using HOSVD to compress the reference solution to multilinear rank~$r$. Observe that the Tucker PRK2-SRRQR method with $h = 10^{-3}$ attains an error comparable error to the quasi-best relative error, within a factor of around 10.

To study how the oblique projector differs from the orthogonal projector, we take the reference solution and compress the tensor at time $t$ using the HOSVD. For the resulting tensor $ \in \mathcal M_rY$, we then compute the relative projection errors of the vector field, that is,  $\|\mathcal{P}_Y^\angle[F(Y)] - F(Y)\|_F / \|F\|_F$ and $\|\mathcal{P}_Y [F(Y)] - F(Y)\|_F / \|F\|_F$. The right panel of Figure~\ref{fig:ac-3d} shows that these errors differ by a constant factor, in line with the statement of Lemma~\ref{lemma: Tucker projection error}.

Finally, Table~\ref{table: performance PRK 3D_AC} compares (Tucker) PRK2-SRRQR and PRK2-ARP with PRK2 (using an orthogonal projector). Due to the rank increase of the nonlinearity (the multilinear rank of $F(Y)$ is at most $2r + r^3$ when $Y$ has multilinear rank-$r$), PRK2 quickly loses its appeal as $r$ increases. On the other hand, the execution times of PRK2-SRRQR and PRK2-ARP grow much more slowly with $r$ and attain a significant advantage already for $r= 5$. Interestingly, PRK2-ARP attains a notable smaller error
than PRK2-SRRQR, which remains very close to the error of PRK.
 \begin{table}[h]
\centering
\begin{tabular}{l rrrr} 
\toprule
Method & Rank=2 & Rank=3 & Rank=4& Rank=5 \\
\midrule
\multicolumn{4}{l}{\textit{Relative errors at final time}} \\ \midrule

PRK2         & 0.2891 & 0.2372 &0.0383&0.0296 \\ 
PRK2-SRRQR &0.3062 &0.2568&0.0654 & 0.0447 \\ 
PRK2-ARP &0.2891& 0.2374 & 0.0391 & 0.0301 \\ 
\midrule 
\multicolumn{4}{l}{\textit{Execution time (in seconds)}} \\ \midrule 
PRK2 & 69.26& 98.28 & 145.66 & 368.01 \\ 
PRK2-SRRQR  &86.34 & 92.83& 111.31 & 129.04 \\ 
PRK2-ARP  & 88.36&96.94 &120.79 & 133.59\\ 
\bottomrule
\end{tabular}
\caption{Performance comparison of Tucker PRK / PRK-DEIM methods applied to the 3D Allen--Cahn equations with step size $h=10^{-3}$ and final time $T=1$.}
\label{table: performance PRK 3D_AC}
\end{table}

\section{Conclusions}

This work has introduced a continuous-time framework, DLRA-DEIM, that integrates data-sparse oblique projections into dynamical low-rank approximation. By interpreting the resulting discontinuous dynamics as a Filippov differential inclusion, we established existence, exactness, and error bounds analogous to those of standard DLRA. For the specific case of QDEIM, the inclusion admits an explicit convex-polytope characterization. On the algorithmic level, we proposed PRK-DEIM and PERK-DEIM, a new family of projected (exponential) Runge–Kutta methods combined with DEIM-based projections. These methods retain the convergence order of existing orthogonally projected (exponential) Runge–Kutta schemes while offering substantial computational savings. The approach extends naturally to low-order tensor differential equations in Tucker format, and its generalization to tensor-trains or general tensor networks presents a promising direction for future work.

 \bibliographystyle{plain}
   \bibliography{Bibliography}

 \appendix
\textbf{}

\section{Proof of Lemma~\ref{lemma:suiu}}
\label{Appendix:Q-DEIM without tie-breaking rule}

\emph{Proof of $\mathcal{I}_U\subset \mathcal{S}_U$.}
Let $S\in \mathcal{I}_U$. By suitably reordering the rows of $U$, we may assume without loss of generality that
$S = [\mathbf{e}_1,\ldots,\mathbf{e}_r]$, that is, the row indices $p_1 = 1,\ldots, p_r  =r$ are a feasible selection of Q-DEIM without tie-breaking. To establish $S \in \mathcal{S}_U$, we will show the following: For every
$\epsilon>0$, there exists an orthonormal matrix $W$ such that $\|W-U\|_F \leq C\epsilon$ for some constant $C$ and Q-DEIM \emph{with} tie-breaking (Algorithm~\ref{alg:QDEIM}) applied to $W$ returns $S$. In fact, we will show that the row of maximal norm is uniquely determined in each iteration of Algorithm~\ref{alg:QDEIM} applied to $W$ and, hence, no tie breaking is needed.

    Because QDEIM is invariant under the choice of orthonormal basis, by taking $U^\top$ as the $R$ factor in the QR decomposition of $U^\top$, we may assume that
    \begin{equation}
    \label{Matrix U}
        U:=[u_1,u_2,\ldots, u_r]=\begin{bmatrix}
\alpha_1 & 0       & 0       & \cdots & 0       \\
*        & \alpha_2 & 0       & \cdots & 0       \\
*        & *       & \ddots  &        & \vdots  \\
\vdots   & \vdots  &         & \alpha_{r-1} & 0 \\
*        & *       & \cdots  & *      & \alpha_r\\
\vdots   & \vdots  &         & \vdots & \vdots \\
*        & *       & *   & *      & * \\
\end{bmatrix}.
    \end{equation} By the orthonormality of $U$ and the fact that $p_i = i$ for $i = 1,\ldots,r$ is a feasible choice of indices, we have that \[1\geq \alpha_1\geq \alpha_2\geq \cdots \geq \alpha_r> 0\] and the row $[\alpha_i,0, \ldots,0]\in \R^{r-i+1}$  is a row of largest norm within the submatrix $[u_{i},\ldots,u_r]\in \R^{(n)\times (r-i+1)}$ for $i=1,\ldots,r$.  
    
We first address the case that $1> \alpha_1$. We will construct an orthonormal matrix $W_0=[w_1,w_2,\ldots ,w_r]$ that is $O(\epsilon)$-close to $U$, has the same lower triangular shape~\eqref{Matrix U}, but the diagonal is guaranteed to be strictly decreasing ($1> \tilde{\alpha}_1>  \tilde{\alpha}_2> \cdots > \tilde{\alpha}_r> 0$), and $[\tilde{\alpha}_i,0, \ldots,0]\in \R^{r-i+1}$ has \emph{the} largest row norm in the matrix $[w_{i},\ldots,w_r]\in \R^{n\times (r-i+1)}$ for $i=1,\ldots,r$. One easily verifies  that in every step $i$ of Algorithm $\ref{alg:QDEIM}$ applied to such a matrix $W$ the row of maximal norm is uniquely determined to be the $i$th row and, in particular, the algorithm returns $S$.

To construct $W$, we perturb each column of $U$ from the last to the first. After $i$ steps of our procedure, the perturbed matrix
$W_{r-i}$ takes for any sufficiently small $\epsilon > 0$ the form  
\begin{equation} \label{eq:Wrmi}
 W_{r-i}=[\tilde{u}_1,\ldots,\tilde{u}_{r-i},w_{r-i+1},w_{r-i+2},\ldots, w_r]= \begin{bmatrix}
\hat{\alpha}_1 & 0       & 0       & \cdots & 0    &0   \\
*        & \ddots & 0       & \cdots & 0     &0  \\
*        & *       & \hat{\alpha}_{r-i} &        & \vdots &0 \\
\vdots   & \vdots  &         & \tilde{\alpha}_{r-i+1} & 0 &0\\
*        & *       & \cdots  & *      & \ddots&0\\
\vdots   & \vdots  &         & \vdots & \ddots &\tilde{\alpha}_{r}\\
*        & *       & *   & *      & * &*\\
*        & *       & *   & *      & * &*\\
\end{bmatrix}
\end{equation}
and it has the following properties:
\begin{itemize}
 \item $W_{r-i}$ is orthonormal;
    \item there exists a constant $C \ge 0$ such that $\|w_i-u_i\|_2\leq C\epsilon$ and $\|\Tilde{u}_i-u_i\|_2\leq C\epsilon$;
    \item $1>\tilde{\alpha}_{r-i+1}>\tilde{\alpha}_{r-i+2}>\ldots >\tilde{\alpha}_r>0$;
    \item $[\tilde{\alpha}_j,0, \ldots,0]\in \R^{r-j+1}$ is \emph{the} row of largest norm in $[w_{j},\ldots,w_r]\in \R^{n\times (r-j+1)}$ for $j=r-i+1,\ldots,r$.
\end{itemize}
Here, and in the following, $C$ denotes a generic, positive constant (independent of $\epsilon$).

Trivially, $W_r = U$ satisfies the properties above for $i = 0$.
We now explain how to perform step $i+1$ for $i\ge 0$ starting from the matrix $W_{r-i}$ with the properties listed above.
First, we construct $w_{r-i}$ by setting
$$w_{r-i}:=\frac{\tilde{u}_{r-i}+C_{r-i}\epsilon\cdot \pmb e_{r-i}}{
               \|\tilde{u}_{r-i}+C_{r-i}\epsilon\cdot \pmb e_{r-i}\|_2}=:[\underbrace{0,\ldots,0}_{\R^{r+i-1}},\tilde{\alpha}_{r-i}, \underbrace{*,\ldots,*}_{\in \R^{n-r-i}}]^\top.$$
Because $\hat \alpha_{r-i}$, $\tilde \alpha_{r-i+1}$ are $O(\epsilon)$-close to
$\alpha_{r-i}$, $\alpha_{r-i+1}$ and 
$1> \alpha_{r-i} \ge \alpha_{r-i+1}$, we may assume that $1> \hat \alpha_{r-i}$ (by choosing $\epsilon$ sufficiently small) and, thus, we can find $C_{r-i} = O(1)$ such that
$$1>\tilde{\alpha}_{r-i}>\tilde{\alpha}_{r-i+1}.$$
Moreover, recall that $[\alpha_{r-i},0, \ldots,0]\in \R^{i+1}$ was a row of largest norm in the corresponding submatrix of $U$. Choosing $C_{r-i}$ sufficiently large thus ensures at the same time that $[\tilde{\alpha}_{r-i},0, \ldots,0]$ is \emph{the} row of largest norm in the perturbed submatrix. Because the last $n-r-i$ elements of $w_{r-i}$ and $\tilde{u}_{r-i}$ are equal up to scaling, the vector $w_{r-i}$ remains orthogonal to $w_{r-i+1},\ldots, w_r$. In contrast, the orthogonality to $\tilde u_1,\ldots, \tilde u_{r-i-1}$ is lost and needs to be restored by the Gram--Schmidt process:
\[
\hat{u}_j=\frac{\tilde{u}_j-\langle \tilde{u}_j, w_{r-i}\rangle w_{r-i}}{\|\tilde{u}_j-\langle \tilde{u}_j,w_{r-i}\rangle w_{r-i}\|_2}, \quad j = 1,\ldots, r-i-1.
\]
Because $w_{r-i}$ is $O(\epsilon)$-close to $\tilde u_{r-i}$, which is orthogonal to $\tilde{u}_j$, it follows that
$\hat{u}_j$ remains $O(\epsilon)$-close to $\tilde{u}_j$ and, thus, $O(\epsilon)$-close to ${u}_j$. In summary, the matrix
$$W_{r-i-1}=[\hat{u}_1,\ldots,\hat{u}_{r-i-1},w_{r-i},w_{r-i+1},\ldots, w_1]$$
has the form~\eqref{eq:Wrmi} and satisfies all the properties above, with $i$ replaced by $i+1$. Repeating the described construction from $i = 0$ to $i = r-1$ gives the desired matrix $W_0$.

It remains to 
discuss the case that some diagonal elements in~\eqref{Matrix U} are equal to one. Let $j \ge 1$ be such that $1 = \alpha_1 = \cdots = \alpha_j$ with either $j= r$ or $\alpha_j > \alpha_{j+1}$. Note that this implies that the first $j$ columns of $U$ are the first $j$ unit vectors. Setting $i = r-j$, we can still perform the first $i$ steps of the procedure as described above, leading to
$$W_{r-i}= \left[ \begin{array}{ccc|cccccc}
                   1 & 0       & 0       & \cdots & 0    &0   \\
0        & \ddots & 0       & \cdots & 0     &0  \\
0        & 0       & 1 &        & \vdots &0 \\ \hline
\vdots   & \vdots  &  0       & \tilde{\alpha}_{r-i+1} & 0 &0\\
0        & 0       & 0 & *      & \ddots &0\\
\vdots   & \vdots  & \vdots        & \vdots & \ddots &\tilde{\alpha}_{r}\\
0        & 0       & 0   & *      & * &*\\
0        & 0       & 0   & *      & * &*\\
                  \end{array}  \right] =: \left[ \begin{array}{c|c} I_{r-i} & 0 \\ \hline 0 & \check W_{r-i}  \end{array} \right].$$ 
The condition $r < n$ implies that there is a nonzero vector $\beta_{r-i} \in \R^{n-r+i}$ that is orthogonal to the columns of the $(n-r+i) \times i$ matrix $\check W_{r-i}$. We now obtain $W_{r-i-1}$ by replacing column $r-i$ of $W_{r-i}$ with
\[
 w_{r-i} = \big[\underbrace{0,\ldots,0}_{\R^{r+i-1}},\tilde{\alpha}_{r-i}, \underbrace{\tilde \beta^\top_{r-i}}_{\in \R^{n-i-r}}\big]^\top,
\]
where $\tilde \beta_{r-i}$ is a nonzero scalar multiple of $\beta_{r-i}$ and $0<\tilde{\alpha}_{r-i}<1$ is such that $w_{r-i}$ has norm one. By construction, $W_{r-i}$ is orthonormal. By choosing $\|\tilde \beta_{r-i}\|_2= O(\epsilon)$ sufficiently small, one can ensure that $\tilde{\alpha}_{r-i} > \tilde{\alpha}_{r-i+1}$ and that $[\tilde{\alpha}_{r-i},0,\ldots,0]$ remains the norm of largest row in the respective submatrix. In other words, $W_{r-i-1}$ satisfies all the properties listed above. By continuing this procedure, we again obtain the desired matrix $W_0$.

\emph{Proof of $\mathcal{S}_U\subset \mathcal{I}_U$.}
Let $S\in \mathcal{S}_U$.  By the definition~\eqref{eq:vicinity} of $\mathcal{S}_U$, there exists a sequence of orthonormal matrices $U_i \to U$ such that the corresponding selection matrices $S^{\mathsf{Q}}_{U_i}$ returned by Q-DEIM satisfy $S^{\mathsf{Q}}_{U_i} \to S$.
Because the set of selection matrices is finite, we may assume $S^{\mathsf{Q}}_{U_i}=S$ for all $i$ by neglecting the first terms of the sequence.
Let $p_1,\ldots,p_r$ denote the indices corresponding to $S$. Then
$$\|U_i^{k-1}(p_k,:)\|_2\geq\|U_i^{k-1}(j,:)\|_2 \quad \forall j = 1,\ldots,n, \quad k = 1,\ldots,r.$$ Using the continuity of Q-DEIM for a \emph{fixed} selection of rows (see proof of Theorem~\ref{thm:QDEIM}), we conclude that $$\|U^{k-1}(p_k,:)\|_2=\lim_{i\rightarrow \infty}\|U_i^{k-1}(p_k,:)\| \geq\lim_{i\rightarrow \infty}\|U_i^{k-1}(j,:)\|_2=\|U^{k-1}(j,:)\|.$$ Thus, $p_1,\ldots,p_r$ satisfy~\eqref{eq:rule} and $S\in \mathcal{I}_U$.

\section{Proof of Lemma~\ref{lemma: limit_f}}
\label{App: supplementary lemma}

In the following, we provide details on the proof of
Lemma~\ref{lemma: limit_f}, closely following the proof of~\cite[Lemma 1]{Paden1987calculus}. For this purpose, we first state two elementary facts.
\begin{lemma}[Lemma A.1 in \cite{Paden1987calculus}]
    \label{lem:decreasing-compact-sets}
    Let $\{E_j\}_{j \in \mathbb{N}}$ be a nested sequence of compact subsets of $\Rmn$ such that 
    $
    E_{j+1} \subseteq E_j
    $ for every $j \in \mathbb{N}$.
    Then
    $$\bigcap_{j \in \mathbb{N}} \textup{conv}(E_j)=\textup{conv}\!\Bigl(\,\bigcap_{j \in \mathbb{N}} E_j\Bigr).$$
    \end{lemma}
\begin{lemma}[Equation A11 in \cite{Paden1987calculus}] \label{lemma:aux2} For $f:\M_r \rightarrow \Rmn$ it holds that 
    \begin{align*}&\bigcap_{k\in \mathbb{N}}\Big\{X\in \Rmn \mid \exists \{Y_i\}\subseteq (Y+\tfrac{1}{k}\mathcal{B})\cap \mathcal{M}_r\text{ s.t. } X=f(Y_i)\rightarrow X\Big\}\\&=\{X\in \mathbb{R}^{m\times n} \mid \exists\{Y_i\}\subseteq \M_r \text{ s.t. } Y_i\rightarrow Y \text{ and } f(Y_i)\rightarrow X\}.\end{align*}
\end{lemma}

\begin{proof}[Proof of Lemma \ref{lemma: limit_f}]
The statement of the lemma follows from the chain of equalities 
\begin{align*}
    \mathcal{F}(Y)&=\bigcap_{\epsilon>0}\overline{\text{conv}}\big[f\big((Y+\epsilon\mathcal{B})\cap \mathcal{M}_r)\big)\big] =\bigcap_{k\in \mathbb{N}}\overline{\text{conv}}\big[f\big((Y+\tfrac{1}{k}\mathcal{B})\cap \mathcal{M}_r\big)\big]
    \\&=\bigcap_{k\in \mathbb{N}}{\text{conv}}\big[\overline{f\big((Y+\tfrac{1}{k}\mathcal{B})\cap \mathcal{M}_r\big)}\big]
     ={\text{conv}}\bigcap_{k\in \mathbb{N}}\overline{f\big((Y+\tfrac{1}{k}\mathcal{B})\cap \mathcal{M}_r\big)}
     \\&={\text{conv}}\bigcap_{k\in \mathbb{N}}\Big\{X\in \Rmn \mid \exists Y_i\subseteq (Y+\tfrac{1}{k}\mathcal{B})\cap \mathcal{M}_r\text{ s.t. } f(Y_i)\rightarrow X\Big\}
     \\&=\textup{conv}\{X\in \mathbb{R}^{m\times n} \mid \exists\{Y_i\}\subseteq \M_r \text{ s.t. } Y_i\rightarrow Y \text{ and } f(Y_i)\rightarrow X\}.
\end{align*}The third equality follows from Assumption \ref{ass: bounded F}, implying that $f\big((Y+\epsilon\mathcal{B})\cap \mathcal{M}_r)\big)$ is bounded.
The fourth equality follows from Lemma~\ref{lem:decreasing-compact-sets} because 
 $$ \overline{f\big((Y+\tfrac{1}{k+1}\mathcal{B})\cap \mathcal{M}_r)\big)}\subseteq \overline{f\big((Y+\tfrac{1}{k}\mathcal{B})\cap \mathcal{M}_r)\big)},$$  and $\overline{f\big((Y+\tfrac{1}{k+1}\mathcal{B})\cap \mathcal{M}_r)\big)}$ is a compact subset of $\mathbb{R}^{m\times n}$ for all $k$. The last  equality is the statement of Lemma~\ref{lemma:aux2}.
\end{proof}

\section{Existence of solutions to differential inclusions}
\label{App: DI}

In the following, we recall a result from~\cite{Haddad1981Monotone} on the existence of a solution to a differential inclusion $\dot{X}(t) \in G(X(t))$ for a
set-valued map $G:\M_r\rightarrow \mathcal{P}(\R^{m \times n})$. For this purpose, we require that $G$ is upper semicontinuous, that is, the following statement holds for every $Y\in \M_r$: For any open set $N$ containing $G(Y)$ there exists a neighborhood $M$ of $Y$ such that $G(M)\subseteq N$; see~\cite[Definition 1]{Aubin1984Differential}. In the special case when $G$ has convex and compact values, an equivalent definition is that the following statement holds for every $Y\in \M_r$: If $\{Y_i\}\subseteq \M_r$, $ Y_i \to Y $ and $ G(Y_i) \ni Z_i \to Z $ then $ Z \in G(Y)$ \cite[page 3714]{Ledyaev2007Nonsmooth}; see~\cite[Lemma 5.15]{teel2012hybrid} for a proof.
\begin{lemma}
    \label{lemma: F_usc}
    Suppose that Assumption~\ref{ass: bounded F} holds. Then the set-valued map $\mathcal{F}$ defined in~\eqref{eq: filippov inclusion} is upper semicontinuous.
\end{lemma}
\begin{proof}
Using Assumption~\ref{ass: bounded F}, the definition~\eqref{eq: filippov inclusion} readily implies that $\mathcal{F}$ has convex and compact values. Now, for $Y\in \M_r$ consider any sequences $\{Y_i\}\subseteq \M_r$, $\{Z_i\} \subseteq \R^{m\times n}$ such that $ Y_i \to Y$ and $ \mathcal{F}(Y_i) \ni Z_i \to Z $. We aim to establish $Z\in \mathcal{F}(Y)$. Choose any $\epsilon>0$. Because
   $Y_i\rightarrow Y$, the inclusion
    $$ (Y_i+\tfrac{\epsilon}{2}\mathcal{B})\cap\M_r\subseteq (Y+\epsilon \mathcal{B})\cap\M_r$$ holds for sufficiently large $i$. Therefore,
    $$Z_i\in\mathcal{F}(Y_i)=\bigcap_{\hat{\epsilon}>0}\overline{\text{conv}} \big[f((Y_i+\hat{\epsilon} \mathcal{B})\cap \M_r)\big] \subseteq\overline{\text{conv}}\big[f((Y_i+\tfrac{\epsilon}{2}\mathcal{B})\cap \M_r)\big]\subseteq \overline{\text{conv}}\big[f((Y+\epsilon \mathcal{B})\cap \M_r)\big].$$ 
    By taking $i\rightarrow \infty$ and because of the closure, we get $Z\in  \overline{\text{conv}}\big[f\big((Y+\epsilon \mathcal{B})\cap \M_r\big)\big].$ Since this hold for arbitrary $\epsilon>0$, the definition of $\mathcal{F}(Y)$ implies that $Z\in\mathcal{F}(Y)$. 
    Therefore, $\mathcal{F}$ is upper semicontinuous.
\end{proof}

Now we are ready to state a general result, \cite[Corollary 1.1]{Haddad1981Monotone} with \cite[Remark 3]{Haddad1981Monotone}, on differential inclusions. This theorem is originally stated for a general locally compact subset $\mathcal{X}$ using Bouligand contingent cones. As a smooth, finite-dimensional manifold, $\mathcal{X}=\M_r$ is locally compact and Bouligand contingent cones coincide with tangent spaces, see~\cite[Ch. 7]{Aubin1984}.

\begin{theorem}
\label{thm: existant}
 Let $G: \M_r \rightarrow \mathcal{P}(\mathbb{R}^{m\times n})$ be upper semicontinuous such that every set returned by $G$ is non-empty, compact and convex. Then, for every $X_0 \in \M_r $ there exists $t_{\max} > 0$ and a Lipschitz continuous function $X: [0,t_{\max}) \to \M_r $ such that $X(0) = X_0$
and $\dot{X}(t) \in G(X(t))$ almost everywhere in $[0,t_{\max})$ if and only if the condition $G(X) \cap T_X\M_r  \neq \emptyset$ holds for all $X \in \M_r$. Furthermore, for any $R > 0$ such that
$$
U_R := \{ X \in \M_r  \mid \|X - X_0\|_F \leq R \}
$$is closed and $G(Y)$ is bounded for every $Y\in U_R$ by $M\geq 0$, one can choose $t_{\max} > R/M$.
\end{theorem}

\section{Proof of Lemma \ref{lemma: projected stages error}}

\label{App: proof of projected stages error}
\begin{proof}[Proof of Lemma \ref{lemma: projected stages error}]
The proof is a direct extension of the one for~\cite[Lemma 7]{kieri2019projection}. In the following, we provide details of this extension, mostly for completeness and to keep track of the modified constants. 

Starting with $Z_1=Y_i=\tilde{Z}_1$, we have $\Tr(Z_1) = Z_1$ and, therefore,
$$\|Z_1 - \tilde{Z}_1\|_F = 0, \quad \big\|\mathcal P_{\Tr(Z_1)}^{\angle}{F(\Tr(Z_1))} - F(\tilde Z_1)\big\|_F \leq C_{m,r} C_{n,r} \cdot \varepsilon_r.$$
The rest of the proof proceeds by induction. Assume that the statement of the lemma holds until stage $j-1$. The constants used depend on that index $j-1$, but we omit writing the dependence to avoid notational overhead.
Using the short-hand notations
$$P_{\ell}^{\angle} = \mathcal P^{\angle}_{\mathcal T_r(Z_{\ell})}, \quad F_{\ell} = F(\mathcal T_r (Z_{\ell})), \quad \tilde F_{\ell} = F(\tilde Z_{\ell}),$$
we obtain from definitions \eqref{eq: Runge--Kutta method} and \eqref{eq: PRK-DEIM}:
$$\|Z_j - \tilde Z_j\|_F \leq h \sum_{\ell = 1}^{j-1} | a_{j \ell} | \|P_{\ell}^{\angle} F_{\ell} - \tilde F_{\ell}\|_F.$$
By the induction hypothesis, it follows that
$$\|Z_j - \tilde Z_j\|_F \leq h \cdot c_2 \cdot \max_{\ell}{|a_{j \ell}|} \cdot \left((j-1) \cdot c_4 \cdot C_{m,r} C_{n,r} \cdot \varepsilon_r + h^{\tilde q_2 +1} + \ldots + h^{\tilde q_{j-1} + 1} \right).$$
Noting that $q_2 = \tilde q_2 \leq \tilde q_3 \leq \ldots \leq \tilde q_{j-1}$, we absorb higher powers of $h$ and the above expression simplifies to
$$\|Z_j - \tilde Z_j\|_F \leq c_1 \cdot h \cdot (c_3 \cdot C_{m,r} C_{n,r} \cdot \varepsilon_r + h^{q_2 + 1}),$$
which is the first inequality claimed in the lemma for stage $j$.

To address the second inequality, we start with the triangle inequality
\begin{equation} \label{eq:triangular}
 \|P_j^{\angle} F_j - \tilde F_j\|_F \leq \|P_j^{\angle} F_j - F_j\|_F + \|F_j - \tilde F_j\|_F \leq C_{m,r} C_{n,r} \cdot \varepsilon_r + L \|\Tr(Z_j) - \tilde Z_j\|_F.
 \end{equation}
By \cite[Lemma 4.1]{hackbusch2016new}, the last term satisfies
$$\|\Tr(Z_j) - \tilde Z_j\|_F \leq \|\Tr(\tilde Z_j) - \tilde Z_j\|_F + 2 \|\tilde Z_j - Z_j\|_F.$$
The classical Runge--Kutta stages verify $\tilde Z_j = \phi_F^{c_j h}(Y_i) + h \Delta$ where $\|\Delta\|_F \le C_L h^{q_j}$ and $c_j = \sum_{\ell = 1}^s a_{j \ell}$. 
Therefore,
\begin{align*}
\|\Tr(\tilde Z_j) - \tilde Z_j\|_F &\leq \min_{X \in \bar{\mathcal M}_r} \|X - \phi_F^{c_j h}(Y_i) + h \Delta\|_F \\
&\leq \min_{X \in \bar{\mathcal M}_r} \|X - \phi_F^{c_j h}(Y_i)\|_F + \|h \Delta\|_F \\
&\leq \varepsilon_r h \max \{1, e^{\ell c_j h}\} + C_L h^{q_j+1},
\end{align*} 
where we used the error bound \eqref{eq: error made by DLRA} in the last inequality.
Inserting these bounds into~\eqref{eq:triangular} gives
\begin{align*}
    \|P_j^{\angle} F_j - \tilde F_j\|_F &\leq C_{m,r} C_{n,r} \cdot \varepsilon_r + L h \big( 2 c_1 (c_3 \cdot C_{m,r} C_{n,r} \cdot \varepsilon_r + h^{q_2+1}) + \varepsilon_r \max \{1, e^{\ell c_j h} \} + C_L h^{q_j} \big) \\
    &\leq c_2 \left(c_4 \cdot C_{m,r} C_{n,r} \cdot \varepsilon_r + h^{\tilde q_j+1} \right),
\end{align*}
which is the second bound claimed in the lemma for stage $j$.
\end{proof}

\end{document}